\documentclass[reqno]{conm-p-l}

\newtheorem{theorem}{Theorem}[section]

\newtheorem{lemma}[theorem]{Lemma}
\newtheorem{corollary}[theorem]{Corollary}

\newtheorem{definition}[theorem]{Definition}

\newtheorem{remark}[theorem]{Remark}

\numberwithin{equation}{subsection}

\begin{document}

\title[Invertibility of the Spectral Operator and the Riemann Hypothesis]{Fractal Complex Dimensions, Riemann Hypothesis and Invertibility of the Spectral Operator}

\author{Hafedh Herichi }
\address{Department of Mathematics, University of California, Riverside, CA 92521-0135}
\curraddr{}
\email{herichi@math.ucr.edu}
\thanks{}

\author{Michel\,L.\,Lapidus}
\address{Department of Mathematics, University of California, Riverside, CA 92521-0135}
\email{Lapidus@math.ucr.edu}
\thanks{The work of M.\,L.\,Lapidus was partially supported by the US National Science foundation under the research grant DMS-1107750, as well as by the Institut des Hautes Etudes Scientifiques (IHES) where the second author was a visiting professor in the Spring of 2012 while this paper was written.}

%    General info
\date{August the 4th, 2012.}

\dedicatory{}


\begin{abstract}
A spectral reformulation of the Riemann hypothesis was obtained in \cite{LaMa2} by the second author and H.\,Maier in terms of an inverse spectral problem for fractal strings.\,The inverse spectral problem which they studied is related to answering the question \textquotedblleft Can one hear the shape of a fractal drum?\textquotedblright and was shown in \cite{LaMa2} to have a positive answer for fractal strings whose dimension is $c\in(0,1)-\{\frac{1}{2}\}$ if and only if the Riemann hypothesis is true.\,Later on, the spectral operator was introduced semi-heuristically by M.\,L.\,Lapidus and M.\,van Frankenhuijsen in their development of the theory of fractal strings and their complex dimensions \cite{La-vF2, La-vF3} as a map that sends the geometry of a fractal string onto its spectrum.\,In this survey article, we focus on presenting the results obtained by the authors in \cite{HerLa1} about the invertibility (in a suitable sense) of the spectral operator, which turns out to be intimately related to the critical zeroes of the Riemann zeta function.\,More specifically, given any $c\geq0$, we show that the spectral operator $\mathfrak{a}=\mathfrak{a}_{c}$, now precisely defined as an unbounded normal operator acting in an appropriate weighted Hilbert space $\mathbb{H}_{c}$, is `quasi-invertible' (i.e., its truncations are invertible) if and only if the Riemann zeta function $\zeta=\zeta(s)$ does not have any zeroes on the vertical line $Re(s)=c$.\,It follows, in particular, that the associated inverse spectral problem has a positive answer for all possible dimensions $c\in (0,1)$, other than the mid-fractal case when $c=\frac{1}{2}$, if and only if the Riemann hypothesis is true.\,Therefore, in this latter result from \cite{HerLa1}, a spectral reformulation of the Riemann hypothesis is obtained from a rigorous operator theoretic point of view, thereby further extending the earlier work of the second author and H.\,Maier in their study of the inverse spectral problem.
\end{abstract}
\maketitle
\newpage \setcounter{tocdepth}{2}
\tableofcontents
\begin{center}
\section{Introduction}
\end{center}

\hspace*{3mm}In \cite{LaMa2}, a spectral reformulation of the Riemann hypothesis was obtained by M.\,L.\,Lapidus and H. Maier, in terms of a family of inverse spectral problems for fractal strings.\,The inverse spectral problem they studied investigates answering the following question:

\begin{quotation}
\hspace*{1.5cm}\textquotedblleft \emph{Can one hear the shape of a fractal string?}\textquotedblright
\end{quotation}

More specifically,

\begin{quotation}
\textquotedblleft \emph{Let $\mathcal{L}$ be a given standard fractal string whose dimension is $D\in(0,1)$.\,If this string has no oscillations of order $D$ in its spectrum, can one deduce that it is Minkowski measurable $($i.e., that it has no oscillations of order $D$ in its geometry$)$?}\textquotedblright
\end{quotation}
\vspace*{2mm}

\hspace*{3mm}The question turned out to have a positive answer other than in the 'midfractal' case, i.e., for any fractal string whose dimension is $D\in (0,1)-\frac{1}{2}$, if and only if the Riemann hypothesis is true.\,(See \cite{LaMa2}, announced in \cite{LaMa1}.)\,This result provided a resolution for the converse of the modified Weyl--Berry conjecture which was formulated in \cite{La1} and then resolved in the affirmative by M.\,L.\,Lapidus and C.\,Pomerance in \cite{LaPo2} (announced in \cite{LaPo1}) in the case of ordinary fractal strings (i.e., one-dimensional drums with fractal boundary).\,Later on, this work was revisited in the light of the \emph{theory of fractal strings and their complex dimensions} which was developed in \cite{La-vF2, La-vF3} by M.\,L.\,Lapidus and M.\,van\,Frankenhuijsen.\\

\hspace*{3mm}In addition, in \cite{La-vF3}, the \emph{spectral operator} was introduced `semi-heuristic-\\ally' as the map that sends the geometry of a fractal string onto its spectrum.\\

\hspace*{3mm}In our recent joint work \cite{HerLa1}, we provided a precise definition of the spectral operator $\mathfrak{a}$ as well as a rigorous functional analytic framework within which to study its main properties.\,We showed that $\mathfrak{a}=\mathfrak{a}_{c}$ is an unbounded normal operator acting on a suitable scale of Hilbert spaces (roughly, indexed by the Minkowski dimension $c$ in (0,1) of the underlying fractal strings) and precisely determined its spectrum (which turned out to be equal to the closure of the range of values of the Riemann zeta function along the vertical line $Re(s)=c$).\,Furthermore, we introduced a suitable family of truncated spectral operators $\mathfrak{a}^{(T)}$ and deduced that for a given $c\geq0$, the spectral operator $\mathfrak{a}=\mathfrak{a}_{c}$ is quasi-invertible (i.e., each of the truncated spectral operators is invertible) if and only if there are no Riemann zeroes on the vertical line of equation  $Re(s)=c$.\,It follows that the associated inverse spectral problem has a positive answer for all possible dimensions $c\in (0,1)$, other than the mid-fractal case when $c=\frac{1}{2}$, if and only if the Riemann hypothesis is true.\\

\hspace*{3mm}Using, in particular, results concerning the universality of the Riemann zeta function among the class of non-vanishing holomorphic functions, we also showed in \cite{HerLa1} that the spectral operator is invertible for $c>1$, not invertible for $\frac{1}{2}<c<1$, and conditionally (i.e., under the Riemann hypothesis), invertible for $0<c<\frac{1}{2}$.\,Moreover, we proved that the spectrum of the spectral operator is bounded for $c>1$, unbounded for $c=1$, equals the entire complex plane for $\frac{1}{2}<c<1$, and unbounded but, conditionally, not the whole complex plane, for $0<c<\frac{1}{2}$.\,We therefore deduced that four types of (mathematical) phase transitions occur for the spectral operator at the critical values (or \emph{critical fractal dimensions}) $c=\frac{1}{2}$ and $c=1$, concerning the shape of its spectrum, its boundedness (the spectral operator is bounded for $c>1$, unbounded otherwise), its invertibility (with phase transitions at $c=1$ and, conditionally, at $c=\frac{1}{2}$), as well as its quasi-invertibility (with a phase transition at $c=\frac{1}{2}$ if and only if the Riemann hypothesis holds true).\\

\hspace*{3mm}The theory of fractal strings and their complex dimensions investigates \emph{the geometric, spectral and physical properties of fractals} and precisely describes \emph{the oscillations in the geometry and the spectrum of fractal strings}; see, in particular, \cite{La-vF2,La-vF3}.\,Such oscillations are encoded in the complex dimensions of the fractal string, which are defined as the poles of the corresponding geometric zeta function.\,This theory has a variety of applications to number theory, arithmetic geometry, spectral geometry, fractal geometry, dynamical systems, geometric measure theory, mathematical physics and noncommutative geometry; see, for example, \cite{La2, La3, La-vF1, La-vF2, La-vF3, La-vF4, La5}; see, in particular, Chapter 13 of the second edition of \cite{La-vF3} for a survey of some of the recent developments in the theory.\\

\hspace*{3mm}The goal of the present survey article is to give an overview of the spectral reformulation of the Riemann hypothesis obtained in \cite{HerLa1} by studying (from various points of view) the invertibility of the spectral operator $\mathfrak{a}=\mathfrak{a}_{c}$, and to show how this work sheds new light (especially, from an operator theoretic perspective) on the earlier reformulation obtained in \cite{LaMa2} and revisited in \cite{La-vF2, La-vF3}.\\

\hspace*{3mm}In closing this introduction, we note that other aspects of the research memoir (or monograph) \cite{HerLa1} are surveyed in \cite{HerLa3} and \cite{HerLa4}.\,In particular, in \cite{HerLa3}, the emphasis is placed on various kinds of mathematical `phase transitions' in connection with the spectral operator and its spectrum, while in \cite{HerLa4}, the focus is on the `universality' of the spectral operator.\,Finally, in the work in preparation \cite{HerLa2}, we study the operator-valued Euler product representation of the spectral operator, both outside and within the critical strip.\\

\hspace*{3mm}The remainder of this paper is organized as follows.\,In \S2, we briefly review the relevant aspects of the theory of generalized fractal strings and their complex dimensions, along with the corresponding explicit formulas (both in the geometric and spectral settings).\,In \S3, after having discussed the heuristic formulation of the spectral operator provided in \cite{La-vF3}, we precisely define and study the infinitesimal shift $\partial_{c}$ (the differentiation operator in one real variable) in terms of which we in turn define the spectral operator $\mathfrak{a}_{c}$.\,Namely, $\mathfrak{a}_{c}=\zeta(\partial_{c})$ (defined via the measurable functional calculus for unbounded normal operators), where $\zeta=\zeta(s)$ is the classic Riemann zeta function.\,We also determine the spectrum of $\partial_{c}$ and give the explicit representation of the shift group generated by $\partial_{c}$.\,In \S4, we explain in more details the original inverse spectral problem for fractal strings studied in \cite{LaMa1, LaMa2} and state the corresponding results obtained therein.\,We also briefly discuss the associated direct spectral problem for fractal strings studied earlier in \cite{LaPo1, LaPo2}, and place it in the broader context of the (modified) Weyl--Berry conjecture for fractal drums [\textbf{La1--4}].\,In \S5, we introduce the truncated infinitesimal shifts and spectral operators in terms of which we can define two new notions of invertibility of $\mathfrak{a}=\mathfrak{a}_{c}$, namely, quasi-invertibility and almost invertibility.\,After having determined the spectra of the above operators and their truncations, we characterize the quasi-invertibility (as well as the almost invertibility) of $\mathfrak{a}_{c}$.\,In \S6, we use the results of \S5 to deduce the aforementioned spectral reformulation of the Riemann hypothesis (RH) (as well as of almost RH, according to which all but finitely many zeroes of $\zeta(s)$ are located on the vertical line $Re(s)=\frac{1}{2}$).\,Finally, in \S7 (and toward the end of \S6), we mention several open problems and extensions of the above results, as well as very briefly discuss some of the other main results of \cite{HerLa1} (or of \cite{HerLa2}).\\

\hspace*{3mm}At the end of the paper, two appendices are also provided.\,In the first one (Appendix A, \S8) we give an elementary overview of some of the main properties of $\zeta$ and of Riemann's beautiful explicit formula connecting the prime numbers and the zeroes of $\zeta$.\,Moreover, in Appendix B (i.e., \S9), we provide an outline of the proof of two key preliminary results from \cite{HerLa1} (discussed in \S3 of the present paper), namely, the normality of the infinitesimal shift $\partial_{c}$ and the characterization of its spectrum $\sigma(\partial_{c})$:
\begin{equation}
\sigma(\partial_{c})=\{s\in \mathbb{C}:\,Re(s)=c\}.\notag
\end{equation}    

\begin{center}
\section{Generalized Fractal Strings and Their Complex Dimensions} 
\end{center}
\subsection{The geometry and spectra of ordinary fractal strings.}
\hspace*{3mm}In fractal geometry, an \emph{ordinary fractal string} is a bounded open subset $\Omega$ of the real line.\,Such a set is a disjoint union of open intervals, the lengths of which form a sequence
 
\begin{equation}
\mathcal{L}=l_{1},l_{2},l_{3},...
\end{equation}
\text
which we will assume to be infinite.\,Since $vol_{1}(\Omega)=\sum_{j\geq 1} l_{j}<\infty$ (where $vol_{1}$ is the one-dimensional Lebesgue measure on $\mathbb{R}$), we may assume without loss of generality that $\{l_{j}\}_{j\geq1}$ is nonincreasing and tends to zero as $j\to \infty$.\\

\hspace*{3mm}Important information about the geometry of $\mathcal{L}$ is contained in its \emph{geometric zeta function},

\begin{equation}\label{Eq:gezet}
\zeta_{\mathcal{L}}(s)=\sum_{j=1}^{\infty}{l_{j}}^{s},
\end{equation}
\text
where $Re(s)>D_{\mathcal{L}}$.\,Here, $D_{\mathcal{L}}:=\inf\{\alpha\in \mathbb{R}:\,\sum_{j=1}^{\infty}l_{j}^{\alpha}<\infty\}$ is the \emph{dimension} of $\mathcal{L}$;\footnote{It then follows that $\{s\in \mathbb{C}:\,Re(s)>\mathcal{D}_{\mathcal{L}}\}$ is the largest open half-plane on which the series $\sum_{j=1}^{\infty}l_{j}^{s}$ is (absolutely) convergent.} it is called the \emph{abscissa of convergence} of the Dirichlet series $\sum_{j=1}^{\infty}{l_{j}}^{s}$ and coincides with the fractal (i.e., Minkowski or box) dimension\footnote{For the notion of Minkowski (or box) dimension, see, e.g., \cite{Fa}, \cite{Mat}, \cite{La1} or \cite{La-vF2, La-vF3}.\,See also Definition \ref{Def:Minmeas} and Remark \ref{RK:MDabcv}.} of the boundary of $\Omega$.\,Furthermore, $\zeta_{\mathcal{L}}$ is assumed to have a suitable meromorphic extension to an appropriate domain of the complex plane containing the half-plane $\{Re(s)>D_{\mathcal{L}}\}$.\\

\hspace*{3mm}\emph{The complex dimensions} of an ordinary fractal string $\mathcal{L}$, as introduced by the second author and M.\,van Frankenhuijsen in the theory of fractal strings and their complex dimensions, are defined as the poles of the meromorphic extension of $\zeta_{\mathcal{L}}$.\,Interesting information about the geometric, spectral (i.e., vibrational) and dynamical \emph{oscillations} of a fractal string is encoded in both the real parts and imaginary parts of its complex dimensions (see \cite{La-vF2,La-vF3} for more information about the theory of ordinary fractal strings and their complex dimensions; see also Remark \ref{RK:FSTdirec}).

\begin{remark}\label{Rk:dimpart}
In the theory of complex dimensions, an object is called fractal if its geometric zeta function has at least one complex dimension with positive real part.\,$($See $[$\emph{\textbf{La-vF3}}, \S12.2$]$.$)$\,As a result, as expected, all $($non-trivial$)$ self-similar geometries are `fractal'.\,Furthermore, other geometries, which could not be viewed as being fractal according to earlier definitions $($in \cite{Man}$)$, are shown to be `fractal' in this new sense, as desired; this is the case, for example, of the Cantor curve $($or `devil's staircase'$)$ and of Peano's plane-filling curve.\,Moreover, every arithmetic geometry ought to be `fractal', due to the presence of the $($critical$)$ zeroes of the corresponding arithmetic zeta function $($or $L$-function$).$\footnote{Those zeroes are viewed as the poles of the logarithmic derivative of the $L$-function $($for instance, the Riemann zeta function in the case of the elusive space attached to the rational number field $\mathbb{Q}$ and the Riemann zeroes$)$.} See \cite{La-vF2, La-vF3, La5}.
\end{remark}

\begin{remark}\label{RK:FSTdirec}
\label{Ft:com}The theory of fractal strings originated in \emph{[\textbf{La1--4}]}, \emph{[\textbf{La-Po1--3}]}, \emph{[\textbf{LaMa1--2}]} and in the memoir \emph{\cite{HeLa}}.\,It was pursued in many directions since then, by the second author and his collaborators, while the mathematical theory of complex fractal dimensions developed and matured; see the books \cite{La-vF2}, \cite{La-vF3} and \cite{La5}.\,See, especially, Chapter $13$ of the second revised and enlarged edition of \cite{La-vF3} for an overview of the theory and for a number of relevant references, including \cite{HamLa} for the case of random fractal strings, \cite{LaLu1, LaLu2, LaLu-vF1, LaLu-vF2} for the case of nonarchimedean $($or $p$-adic$)$ fractal strings, \cite{LaLeRo, ElLaMaRo} for the study of multifractal strings, as well as \cite{LaPe} and \cite{LaPeWi} where the beginning of a higher-dimensional theory of complex dimensions of fractals is developed, particularly under suitable assumptions of self-similarity.\,$($See also \cite{LaRaZu} for a significantly more general higher-dimensional theory, potentially applicable to arbitrary fractals.$)$
\end{remark}

\hspace*{3mm}The \emph{Cantor string}, denoted by $CS$, and defined as the complement of the Cantor set in the closed unit interval $[0,1]$, is a standard example of an ordinary fractal string:\\

$CS=(\frac{1}{3},\frac{2}{3})\bigcup(\frac{1}{9},\frac{2}{9})\bigcup(\frac{7}{9},\frac{8}{9})\bigcup(\frac{1}{27},\frac{2}{27})\bigcup(\frac{7}{27},\frac{8}{27})\bigcup(\frac{19}{27},\frac{20}{27})\bigcup(\frac{25}{27},\frac{26}{27})\bigcup ...$\\

\text
Here, each length $l_{j}=3^{-j-1}$, $j\geq0$, is counted with a multiplicity $w_{j}=2^{j}$.\,Thus, the geometric zeta function associated to such a string is 

\begin{equation}\label{Eq:zetCS}
\zeta_{CS}(s)=\sum_{j=0}^{\infty}2^{j}.3^{-(j+1)s}=\frac{3^{-s}}{1-2.3^{-s}}\\
\end{equation}
\text
whose set of poles is the set of complex numbers
\begin{equation}\label{Eq:DCS}
\mathcal{D}_{CS}=\{D+in\textbf{p}:\,n\in\mathbb{Z}\},\\
\end{equation}
\text 
where $D=\log_{3}2$ is the dimension of the $CS$ and \textbf{p}=$\frac{2\pi}{\log 3}$.\,This set is called \emph{the set of complex dimensions} of $CS$.\,Note that the real part of these complex numbers is the Minkowski dimension of $CS$ and that the imaginary parts correspond to the oscillatory period \textbf{p} in the volume of the inner tubular neighborhoods of $CS$, as we now explain.\\

\hspace*{3mm}For a given $\epsilon>0$, the \emph{volume of the inner $\epsilon$-tubular neighborhood} of the boundary, $\partial\Omega$, of a fractal string $\mathcal{L}$ is
\begin{equation}
V_{\mathcal{L}}(\epsilon)=vol_{1}\{x\in\Omega:\,d(x,\partial\Omega)<\epsilon\},
\end{equation}
where $vol_{1}$ is the one-dimensional Lebesgue measure on $\mathbb{R}$, as before, and $d(x, \partial\Omega)$ denotes the distance from a point $x\in \mathbb{R}$ to the boundary of $\Omega$.\,In the case of the Cantor string $CS$ and as is shown in [\textbf{La-vF3}, \S1.1.2], we have 
\begin{equation}
V_{CS}(\epsilon)=\frac{2^{-D}\epsilon^{1-D}}{D(1-D)\log 3}+\frac{1}{\log 3}\sum_{n=1}^{\infty}Re\bigg(\frac{(2\epsilon)^{1-D-in\textbf{p}}}{(D+in\textbf{p})(1-D-in\textbf{p})}\bigg)-2\epsilon.\label{Eq:Etb}
\end{equation}
In this manner, the geometric oscillations that are intrinsic to the Cantor set (viewed as the fractal boundary of the Cantor string) are expressed in terms of the underlying complex dimensions.\\

\hspace*{3mm}More generally, another representation of the volume $V_{\mathcal{L}}(\epsilon)$ of the inner tubular neighborhood of a fractal string $\mathcal{L}$ was obtained by using the explicit formulas from [\textbf{La-vF3}, Ch.\,5] (to be presented and discussed later on in this paper, see Theorem \ref{Thm:2}).\,More specifically, under some mild assumptions, the following `fractal tube formula' is established in [\textbf{La-vF3}, Ch.\,8], enabling one to express $V_{\mathcal{L}}(\epsilon)$ as a sum over the complex dimensions of the fractal string $\mathcal{L}$:
\begin{equation}
V_{\mathcal{L}}(\epsilon)=\sum_{\omega\in \mathcal{D}_{\mathcal{L}}(\mathcal{W})}res\bigg(\frac{\zeta_{\mathcal{L}}(s)(2\epsilon)^{1-s}}{s(1-s)};\omega\bigg)+\{2\epsilon\zeta_{\mathcal{L}}(0)\}.\label{Eq:tbform}
\end{equation}
Here, the term in braces is included only if $0\in \mathcal{W}-\mathcal{D}_{\mathcal{L}}(\mathcal{\mathcal{W}})$.\,Furthermore, $\mathcal{D}_{\mathcal{L}}(\mathcal{W})$ denotes the set of \emph{visible complex dimensions} relative to a suitable `window' $\mathcal{W}\subset\mathbb{C}$ (i.e., the set of poles in $\mathcal{W}$ of the meromorphic continuation of $\zeta_{\mathcal{L}}$ to a connected open neighborhood of $\mathcal{W}$); see [\textbf{La-vF3},\,\S1.2.1].\,Moreover, in general, the tube formula (\ref{Eq:tbform}) also contains an error term which can be explicitly estimated as $\epsilon\to 0^{+}$.\\

\hspace*{3mm}If we assume, for the simplicity of exposition, that $0\notin\mathcal{W}$, $1\notin \mathcal{D}_{\mathcal{L}}(\mathcal{W})$, and that all of the visible complex dimensions are simple (i.e., are simple poles of $\zeta_{\mathcal{L}}$), then (\ref{Eq:tbform}) becomes

\begin{equation}\label{Eq:tbf}
V_{\mathcal{L}}(\epsilon)=\sum_{\omega\in \mathcal{D}_{\mathcal{L}}(\mathcal{W})}res\big(\zeta_{\mathcal{L}}(s);\omega\big)\frac{(2\epsilon)^{1-\omega}}{\omega(1-\omega)},
\end{equation}
which is often referred to as a `fractal tube formula'. (See Theorem 8.1 and Corollary 8.3 in \cite{La-vF3}.)\\

\hspace*{3mm}Note that in the above case of the Cantor string $CS$, we have $\mathcal{W}=\mathbb{C}$, $\mathcal{D}_{CS}(\mathcal{W})=\mathcal{D}_{CS}$, and the error term vanishes identically.\,In addition, the resulting exact (or fractal) tube formula (\ref{Eq:Etb}) holds pointwise (rather than just distributionally), in agreement with the pointwise tube formulas also obtained in [\textbf{La-vF3}, \S8.1.1 \& \S8.4].\,Finally, observe that (\ref{Eq:Etb}) follows from (\ref{Eq:tbf}) since (in light of Equations (\ref{Eq:zetCS}) and (\ref{Eq:DCS})) the complex dimensions of $CS$ are simple and have the same residue, $\frac{1}{\log3}$.\\

\hspace*{3mm}We will see shortly that the explicit distributional formulas play an important role in motivating the definition of the spectral operator $\mathfrak{a}_{c}$.\\

\hspace*{3mm}Spectral information (representing the  \emph{frequencies} of the `\emph{vibrations}' of the fractal string) can also be derived.\,Indeed, one can \emph{listen to the sound} of a given ordinary fractal string $\mathcal{L}=\{l_{j}\}_{j=1}^{\infty}$.\,Here, the positive numbers $l_{j}$ denote the lengths of the connected components (i.e., open intervals) of a bounded open set $\Omega$ of the real line $\mathbb{R}$, with (possibly) fractal boundary $\partial\Omega$.\,In fact, spectral information about $\mathcal{L}$ is encoded by its \emph{spectral zeta function}, which is defined as 

\begin{equation}
\zeta_{\mathcal{\nu}}(s)=\sum_{f} f^{-s},
\end{equation}
\text
where $f=kl_{j}^{-1}$ ($k,j=1,2,...$) are the normalized frequencies of $\mathcal{L}$.\,Up to a trivial normalization factor, these are simply the square roots of the eigenvalues of the Laplacian (or free Hamiltonian) on $\Omega$, with Dirichlet boundary conditions on $\partial\Omega$.\,So that, in particular, the associated eigenfunctions are constrained to have nodes at each of the endpoints of the intervals of which the open set $\Omega$ is composed (see, e.g., [\textbf{La1--5, LaPo1--3, LaMa1--2, HeLa, La-vF2, La-vF3}] for more details).\\

\hspace*{3mm}\emph{The geometry and the spectrum} of  $\mathcal{L}$\space are related via the following formula (see [\textbf{La2--3}],\,[\textbf{LaMa2}],\,[\textbf{La-vF3},\,\S1.3]):
\begin{equation}
\zeta_{\mathcal{\nu}}(s) =\zeta_{\mathcal{L}}(s) .\zeta(s),\label{Eq:gs}
\end{equation}
where $\zeta$ is the Riemann zeta function.\,Here, $\zeta_{\mathcal{L}}$ is the \emph{geometric zeta function} of $\mathcal{L}$, defined by $\zeta_{\mathcal{L}}(s)=\sum_{j=1}^{\infty}l_{j}^{s}$, for $Re(s)>D_{\mathcal{L}}$, the \emph{abscissa of convergence} of the Dirichlet series $\sum_{j=1}^{\infty}l_{j}^{s}$ or \emph{dimension} of $\mathcal{L}$ (which coincides with the Minkowski dimension of $\partial\Omega$, see \cite{La2},\,[\textbf{La-vF3}, \S1.2], along with Definition \ref{Def:Minmeas} and Remark \ref{RK:MDabcv} below).\\

\hspace*{3mm}Equation (\ref{Eq:gs}) plays a key role in \emph{connecting the spectrum of a fractal string $\mathcal{L}$ to its geometry} or conversely (and provided no zero of $\zeta$ coincides with a visible complex dimension of $\mathcal{L}$), in \emph{relating the geometry of a fractal string to its spectrum via the Riemann zeta function}.\\

\hspace*{3mm}In hindsight, this relation helps explain the approach to the \emph{direct spectral problem for fractal strings} adopted in [\textbf{LaPo\,1--2}] and the approach to \emph{inverse spectral problems for fractal strings} used in [\textbf{LaMa\,1--2}].\,We stress, however, that a number of technical difficulties had to be overcome in order to formulate and derive the results obtained in those papers.\,In addition, the notion of complex dimension that was only hidden or heuristic in [\textbf{LaPo\,1--2}], [\textbf{LaMa\,1--2}], [\textbf{La\,1--3}] and [\textbf{HeLa}], was developed rigorously in [\textbf{La-vF2}, \textbf{La-vF3}] (and other papers by the authors of these monographs, beginning with \cite{La-vF1} and several earlier IHES preprints) in part to provide a systematic approach (via explicit formulas generalizing  Riemann's explicit formula discussed in Appendix A) to the results on direct and inverse spectral problems obtained in \emph{loc.\,cit.} (See, for example, [\textbf{La-vF3}, Chs.\,6 \& 9].)\\

\begin{remark}\label{RK:Lj}
Various extensions of the factorization formula $($\ref{Eq:gs}$)$ have since been obtained in \cite{Tep1, Tep2, DerGrVo, LalLa}, in the context of analysis on fractals \cite{Ki} and $($single or multi-variable$)$ complex dynamics, using the decimation method \emph{[\textbf{Ram, RamTo, Sh, FukSh, Sab1--3}]} for the eigenfunctions of Laplacians on certain self-similar fractals.
\end{remark}

\hspace*{3mm}A consequence of a special case of the explicit formulas of [\textbf{La-vF2}, \textbf{La-vF3}] applied to the spectrum and the geometry of a fractal string (in the spirit of formula (\ref{Eq:gs})) is that the Riemann zeta function $($as well as a large class of arithmetic zeta functions and other Dirichlet series$)$ cannot have an infinite vertical arithmetic progression of zeroes.\footnote{In the special case of $\zeta$, this result was already obtained by C.\,Putnam [\textbf{Put1,\,Put2]} in the 1950s, via a completely different proof which does not extend to the general case considered in [\textbf{La-vF1},\,\textbf{La-vF2},\,\textbf{La-vF3}].}\,(See [\textbf{La-vF3},\,Chs.\,10\,\&\,11] for a proof of this result and several of its refinements.)\,For instance, applying the aforementioned explicit formulas to the Cantor String $CS$ and assuming that the Riemann zeta function were to vanish at the complex dimensions $D+in\textbf{p}$ of $CS$, where $n\in \mathbb{Z}-\{0\}$, then one can deduce that $CS$ would have to sound the same as a Minkowski measurable fractal string of the same dimension $D=\log_{3}2$.\\

\hspace*{3mm}A similar conclusion can be obtained by considering generalized Cantor strings (with noninteger multiplicities).\,This conclusion is contradicted by the results of [\textbf{La-vF3}, Ch.\,10] according to which such fractal strings always have geometric oscillations (of leading order) in their geometry, from which one deduces the above theorem about the nonexistence of zeroes in infinite arithmetic progressions.\\

\hspace*{3mm}From now on, we will denote by $\zeta_{\mathcal{L}}$ (respectively, $\zeta_{\nu}$) the meromorphic continuation (when it exists) of the geometric zeta function (respectively, of the spectral zeta function) of a fractal string $\mathcal{L}$.

\subsection{Generalized fractal strings and their explicit formulas.}
\hspace*{3mm}Next, we introduce one of our main objects of investigation, the class of \emph{generalized fractal strings}, and some of the mathematical tools needed to study it.\,(See [\textbf{La-vF3}, Ch.\,4].)\\

\hspace*{3mm}A \emph{generalized fractal string} $\eta$ is defined as a local positive or complex measure on $(0,+\infty)$\,satisfying $|\eta|(0,x_{0})=0$,\,for some $x_{0}>0$.\footnote{In short, a positive (or complex) \emph{local measure} on $(0,+\infty)$ is a locally bounded set-function on $(0,+\infty)$ whose restriction to any bounded Borel subset (or equivalently, bounded subinterval) of $(0,+\infty)$ is a standard positive (or complex) measure.\,See, e.g., [\textbf{La-vF3}, \S4.1].}\,Here, the positive (local) measure $|\eta|$ is the variation of $\eta$.\footnote{For an introduction to measure theory, we refer, e.g., to \cite{Coh,Fo}.\,Recall that when $\eta$ is positive, then $|\eta|=\eta$.}\,A standard example of a generalized fractal string can be obtained as the measure associated to an ordinary fractal string  $\mathcal{L}=\{L_{j}\}_{j=1}^{\infty}$ with multiplicities $w_{j}$.\,(Here, $\{L_{j}\}_{j=1}^{\infty}$ denotes the sequence of \emph{distinct} lengths of $\mathcal{L}$, written in decreasing order and tending to zero as $j\to \infty$.)\,Such a measure is defined as
\begin{equation}
\eta_{\mathcal{L}}=\sum_{j=1}^{\infty}w_{j}\delta_{\{L_{j}^{-1}\}}.\label{Eq:measet}
\end{equation}
Note that $\eta_{\mathcal{L}}$ is a generalized fractal string since $|\eta_{\mathcal{L}}|$ does not have any mass on $(0,L^{-1}_{1})$.\,Here and in the sequel, $\delta_{\{x\}}$ is the Dirac delta measure or the unit point mass concentrated at $x>0$.\\

\begin{remark}
In many important situations where an extension of formula $($\ref{Eq:measet}$)$ is used, one should think of the positive numbers $L_{j}$ $($or their analog $\ell_{j}$ in Equation $($\ref{Eq:gezet}$)$ and the discussion preceding it$)$ as \emph{scales} rather as the lengths associated with some concrete geometric object.\,Furthermore, as we will see next in Remark \ref{RK:multgfrac}, the multiplicities $w_{j}$ need not be integers, in general.
\end{remark}

\begin{remark}\label{RK:multgfrac}
In the case of an ordinary fractal string, $w_{j}$ is always integral for any $j\geq 1$.\,However, in general, this multiplicity $($or weight$)$ is not necessarily integral.\,For instance, the \emph{prime string} 
\begin{equation}
\eta_{\mathcal{B}}=\sum_{m\geq1,p}\big(\log p\big)\delta_{p^{m}},\label{Eq:Prmst}
\end{equation}
where $p\in \mathcal{P}$:= the set of all prime numbers, is an example of a generalized fractal string for which $w_{j}=\log p$ is non-integral.\,It is also the measure associated to the non-ordinary fractal string $\mathcal{L}=\{p_{j}^{-m}\}_{j=1}^{\infty}$ with multiplicities $\log p_{j}$, where $p_{j}$ is the j-th prime number written in increasing order.\,Hence, the use of the word `generalized' is well justified for this class of strings.\\
\end{remark}
\hspace*{3mm} Let $\eta$ be a generalized fractal string.\,Its \emph{dimension} is 
\begin{equation}
D_{\eta}:=\inf\Big\{\sigma\in\mathbb{R}: \int_{0}^{\infty}x^{-\sigma}|\eta|(dx)<\infty\Big\}.\label{Eq:Dimeta}
\end{equation}
The \emph{counting function} of $\eta$ is\footnote{More precisely, in order to obtain accurate \emph{pointwise} formulas, one must set $N_{\eta}(x)=\frac{1}{2}(\eta(0,x]+\eta[0,x))$, much as in the pointwise theory of Fourier series.}
\begin{equation}
N_{\eta}(x):=\int_{0}^{x}\eta(dx)=\eta(0,x).
\end{equation}

The \emph{geometric zeta function} associated to $\eta$ is the Mellin transform of $\eta$.\,It is defined as
\begin{equation} 
\zeta_{\eta}(s):=\int_{0}^{\infty}x^{-s}\eta(dx) \mbox{\quad for $Re(s)>D_{\eta}$},
\end{equation}
where $D_{\eta}$ is the dimension of $\eta$ (and is also called the \emph{abscissa of convergence} of the Dirichlet integral $\int_{0}^{\infty}x^{-s}\eta(dx)$).\,As we did in \S2.1, we assume that $\zeta_{\eta}$ has a meromorphic extension to some suitable (open, connected) neighborhood $\mathcal{W}$ of the half-plane $\{Re(s)>D_{\eta}\}$ (see [\textbf{La-vF3}, \S5.3] for more details on how the window $\mathcal{W}$ is defined) and we define the set $\mathcal{D}_{\mathcal{\eta}}(\mathcal{W})$ of visible \emph{complex dimensions} of $\eta$ by\footnote{Since $\zeta_{\mathcal{L}}$ is assumed to be meromorphic, $\mathcal{D}_{\mathcal{L}}(W)$ is a \emph{discrete} $($and hence, at most countable$)$ subset of $\mathbb{C}$.\,Furthermore, since $\zeta_{\mathcal{L}}$ is holomorphic for $Re(s)>D_{\mathcal{L}}$ $($because by definition of $D_{\mathcal{L}}$, the Dirichlet integral $\int_{0}^{\infty}x^{-s}\eta(dx)$ is absolutely convergent there$)$, all the complex dimensions $\omega$ of $\mathcal{L}$ satisfy $Re(\omega)\leq D_{\mathcal{L}}$.}

\begin{equation}
\mathcal{D}_{\mathcal{\eta}}(\mathcal{W}):=\{\omega\in\mathcal{W}:\zeta_{\eta} \mbox{\quad has a pole at $\omega$}\}. 
\end{equation}
\hspace*{3mm}For example, the geometric zeta function of the prime string $\eta_{\mathcal{B}}$ (defined above in Equation (\ref{Eq:Prmst})) is
\begin{equation}
\zeta_{\eta_{\mathcal{\beta}}}(s)=-\frac{\zeta'(s)}{\zeta(s)}\mbox{\quad for $s\in\mathbb{C}$}.
\end{equation}
Therefore, the complex dimensions of $\eta_{\beta}$ are the zeroes of $\zeta$ (each counted with multiplicity one, along with the single and simple pole of $\zeta$ (located at $s=1$).\,We recall that the trivial zeroes of $\zeta$ occur at the values $s=-2n$, for $n=1,\,2,\,3,\,...$\,The nontrivial (or critical) zeroes of the Riemann zeta function, which are located inside the \emph{critical strip} (i.e., inside the region $0<Re(s)<1$ of the complex plane), are conjectured to lie on the vertical line $Re(s)=\frac{1}{2}$; this celebrated conjecture is known as the Riemann hypothesis.\,(See Appendix A.)\\

\hspace*{3mm}The \emph{spectral measure} $\nu$ associated to $\eta$ is defined by
\begin{equation}
\nu(A)=\sum_{k=1}^{\infty}\eta\left(\frac{A}{k}\right),
\end{equation} 
\text
for any bounded Borel set (or equivalently, interval) $A\subset(0,+\infty)$.\,The geometric zeta function of $\nu$ is then called the \emph{spectral zeta function} associated to $\eta$.\\

\hspace*{3mm}Two important generalized fractal strings (within our framework) are \emph{the harmonic generalized fractal string}
\begin{equation}
\mathfrak{h}=\sum_{k=1}^{\infty}\delta_{\{k\}},
\end{equation}
and \emph{the prime harmonic generalized fractal string} defined for each prime $p\in \mathcal{P}$ as
\begin{equation}
\mathfrak{h}_{p}=\sum_{k=1}^{\infty}\delta_{\{p^{k}\}},
\end{equation}
where, as before, $\delta\{.\}$ is the Dirac delta measure.\,They will play a key role in defining the spectral operator and its operator-valued Euler product (see Equations (\ref{Eq:Spr}) and (\ref{Eq:Pf})).\,These strings are related via the multiplicative convolution operation of measures $\ast$ as follows:
\begin{equation}
\mathfrak{h}=\underset{p\in\mathcal{P}}{\ast}\mathfrak{h}_{p}.
\end{equation}
As a result, we have 
\begin{equation}\label{Eq:Eph}
\zeta_{\mathfrak{h}}(s)=\zeta_{\underset{ p\in\mathcal{P}}{\ast\mathfrak{h_{p}}}}(s)=\zeta(s)=\underset{ p\in\mathcal{P}}{\prod}\frac{1}{1-p^{-s}}=\underset{ p\in\mathcal{P}}{\prod}\zeta_{\mathfrak{h}_{p}}(s),
\end{equation}
for $Re(s)>1$.\\

The \emph{spectral zeta function} associated to $\nu$,\space which as we have seen, is defined as the geometric zeta function of $\nu$,\,is related to $\zeta_{\eta}$\space via the following formula (which is the exact analog of Equation (\ref{Eq:gs})):
\begin{equation} 
\zeta_{\nu}(s)=\zeta_{\eta}(s).\zeta(s),\label{Eq:factor}
\end{equation}
where $\zeta$\space is the Riemann zeta function.\,As is recalled in Appendix A, $\zeta$ is well known to have an \emph{Euler product expansion} given by the formula 
\begin{equation}
\zeta(s)=\prod_{p\in \mathcal{P}}(1-p^{-s})^{-1},\mbox{\quad for\,\,} Re(s)>1,\label{eq:Ep}
\end{equation}
\text
where, as before, p runs over the set $\mathcal{P}$ of all prime numbers.\, Note that this Euler product was reinterpreted differently in Equation (\ref{Eq:Eph}) just above.\\

\hspace*{3mm}Throughout their development of the theory of complex dimensions and inspired by Riemann's explicit formula,\footnote{See Appendix A for a discussion of the analogy between Riemann's original explicit formula and the explicit formulas for generalized fractal strings obtained in \cite{La-vF2,La-vF3}.} which expresses the counting function of the number of primes less than some positive real number in terms of the zeroes of the Riemann zeta function, the second author and M.\,van Frankenhuijsen obtained (and extensively used) \emph{explicit distributional formulas} associated to $\eta$.\,In these formulas,\,the k-th distributional primitive (or anti-derivative) of $\eta$\space is viewed as a \emph{distribution},\,acting on functions in the Schwartz space \cite{Schw} on the half-line $(0,+\infty)$.\,(See [\textbf{La-vF3}, Ch.\,5] for a detailed discussion and a precise statement of these explicit formulas, both in their distributional and pointwise form.)

\begin{theorem}\emph{\cite{La-vF2,La-vF3}} \label{Thm:2}
Let $\eta$\,be a languid generalized fractal string.\footnote{A generalized fractal string $\eta$ is said to be \emph{languid} if its geometric zeta function $\zeta_{\eta}$ satisfies some suitable \emph{polynomial growth conditions}; see [\textbf{La-vF3}, \S5.3].}\\Then, for any $k\in\mathbb{Z}$, its  k-th distributional primitive is given by
\begin{align}
\mathcal{P}_{\eta}^{[k]}(x)=&\sum_{\omega\in \mathcal{D}_{\eta}(\mathcal{W})}res\left(\frac{x^{s+k-1}\zeta_{\eta}(s)}{(s)_{k}};\omega\right)+\frac{1}{(k-1)!}\sum_{\begin{subarray}{l}\hspace*{4.5mm}j=0,\\-j\in \mathcal{W}-\mathcal{D}_{\eta}\end{subarray}}^{k-1}{{k-1}\choose{j}}(-1)^{j}x^{k-1-j}\notag\\
                           &.\zeta_{\eta}(-j)+\mathcal{R}_{\eta}^{[k]}(x),
\end{align}                             
where $\omega$ runs through the set $\mathcal{D}_{\eta}(\mathcal{W})$ of visible complex dimensions of $\eta$ and $\mathcal{R}_{\eta}^{[k]}(x)=\frac{1}{2\pi i}\int_{\mathcal{S}}x^{s+k-1}\zeta_{\eta}(s)\frac{ds}{(s)_{k}}$\,is the error term, which can be suitably estimated as $x \to +\infty$ and which, under appropriate hypotheses, vanishes identically $($thereby yielding \emph{exact} explicit formulas, see \emph{[\textbf{La-vF3}, \S5.3 \& \S5.4]}$)$.\footnote{Here, in general, the binomial coefficients ${{k-1}\choose{j}}$ are defined in terms of the gamma function $\Gamma=\Gamma(s)$.\,Moreover, the Pochhammer symbol is defined by $(s)_{k}=s(s+1)...(s+k-1)$, for $k\geq1$, and $(s)_{k}=\frac{\Gamma(s+k)}{\Gamma(s)}$ for any $k\in \mathbb{Z}$.}
\end{theorem}

\hspace*{3mm}The explicit distributional formula stated in Theorem \ref{Thm:2} provides a representation of $\eta$, in the distributional sense, as a sum over its complex dimensions which encode in their real and imaginary parts important information about the (geometric, spectral or dynamical) oscillations of the underlying fractal object.\,Recall also from our discussion in Appendix A that the original explicit formula was first obtained by Riemann in 1858 as an analytical tool to understand the distribution of the primes.\,It was later extended by von Mangoldt and led in 1896 to the first rigorous proof of the Prime Number Theorem,\,independently by Hadamard and de la Vall\'ee Poussin.\,(See \cite{Edw, Ing, Ivi, Pat, Tit}.)\,In [\textbf{La-vF3}, \S5.5], the interested reader can find a discussion of how to recover the Prime Number Theorem, along with a suitable form of Riemann's original explicit formula and its various number theoretic extensions, from Theorem \ref{Thm:2} (and more general results given in [\textbf{La-vF3}, Ch.\,5]).\\

\hspace*{3mm}Note that Theorem \ref{Thm:2} enables us to obtain, in the distributional sense, useful representations of the k-th primitives of $\eta$, for any $k\in \mathbb{Z}$.\,For instance, if we apply it at level $k=0$, we obtain an explicit representation of $\eta$ which is called the \emph{density of geometric states} formula\,(see [\textbf{La-vF3}, \S6.3.1] and Remark \ref{Rk:Exfor}):\\
\begin{equation}
\eta=\sum_{\,\omega\in\mathcal{D_{\eta}(W)\,}}res(\zeta_{\eta}(s);\omega)x^{\omega-1}.\label{Eq:18}
\end{equation}
We also recall that the spectral measure $\nu=\eta\ast\mathfrak{h}$ is itself a generalized fractal string.\,Thus, when applying the explicit formulas (also at level $k=0$), we obtain an explicit formula for $\nu$ which is similar to the \emph{density of spectral states} (\emph{or density of frequencies}\,formula) in quantum physics (see [\textbf{La-vF3}, \S6.3.1]):\\
\begin{align}
\nu
&=\zeta_{\eta}(1)+\sum_{\,\omega\in\mathcal{D_{\eta}(W)\,}}res(\zeta_{\eta}(s)\zeta(s)x^{s-1};\omega)\notag\\
&=\zeta_{\eta}(1)+\sum_{\omega\in \mathcal{D}_{\eta}(W)}res(\zeta_{\eta}(s);\omega)\zeta(\omega)x^{\omega-1},\label{Eq:19}
\end{align}
\text
where $0\notin \mathcal{W}$, $1\notin \mathcal{D}_{\eta}(\mathcal{W})$ and, as above, $res(\zeta_{\eta}(s);\omega)$ denotes the residue of $\zeta_{\eta}(s)$ as $s=\omega$.

\begin{remark}\label{Rk:Exfor}
Note that the explicit expressions for $\eta$ and $\nu$, stated respectively in Equation $($\ref{Eq:18}$)$ and Equation $($\ref{Eq:19}$)$, are given as a sum over the complex dimensions of $\eta$.\,Here, for clarity, we stated these formulas in the case of simple poles and neglected including the possible error terms.
\end{remark}

\hspace*{3mm}Next, we introduce the \emph{spectral operator}, as defined in \cite{La-vF3} and present some of its fundamental properties, which are rigorously studied in \cite{HerLa1} and further discussed in \cite{HerLa2}.
\begin{center}
\section{The Spectral Operator $\mathfrak{a}_{c}$ and the Infinitesimal Shifts $\partial_{c}$} 
\end{center}
\subsection{A `heuristic' definition of $\mathfrak{a}_{c}$.}
\hspace*{3mm}Following, in particular, the work in [\textbf{La1--3, LaPo1--3, LaMa1--2, HeLa}], relating the \emph{spectrum} of certain classes of fractal strings to their \emph{geometry} has been a subject of significant interest to the authors of \cite{La-vF2, La-vF3} throughout their development of the theory of fractal strings and their complex dimensions.\,Motivated by this fact and also the formula recalled in Equations (\ref{Eq:gs}) and (\ref{Eq:factor}), the \emph{spectral operator} was `heuristically' defined in [\textbf{La-vF3}, \S6.3.2]\footnote{By 'heuristically', we mean that the spectral operator and its operator-valued Euler product (see Equations (\ref{Eq:spop}) and (\ref{Eq:EPr})) were defined in [\textbf{La-vF3}, \S6.3.2] without introducing a proper functional analytic framework enabling one to rigorously study their properties and provide conditions ensuring their invertibility.} as \emph{the operator mapping the density of geometric states $\eta$\space to the density of spectral states $\nu$}:\footnote{This is the level $k=0$ version of the spectral operator, in the sense of Theorem \ref{Thm:2} and of the ensuing discussion.}
\begin{equation}\label{Spoplev1}
\eta\longmapsto\nu
\end{equation}

\
\hspace*{3mm}At level k=1,\,it will be defined on a suitable Hilbert space $\mathbb{H}_{c}$, where $c\geq 0$, as the operator mapping the counting function of $\eta$ to the counting function of $\nu=\eta\ast\mathfrak{h}$ (that is, mapping \emph{the geometric counting function} $N_{\eta}$ onto the \emph{spectral counting function} $N_{\nu}$):
\begin{equation}
N_{\eta}(x)\longmapsto\nu(N_{\eta})(x)=N_{\nu}(x)=\sum_{n=1}^{\infty}N_{\eta}\left(\frac{x}{n}\right).
\end{equation}

\hspace*{3mm}Note that under the change of variable $x=e^{t}$, where $t\in \mathbb{R}$ and $x>0$, one can obtain an \emph{additive} representation of the spectral operator $\mathfrak{a}$,
\begin{equation}\label{Eq:spop}
f(t)\mapsto \mathfrak{a}(f)(t)=\sum_{n=1}^{\infty}f(t-\log n),
\end{equation}
and of its \emph{operator-valued Euler factors} $\mathfrak{a}_{p}$
\begin{equation}
f(t)\mapsto \mathfrak{a}_{p}(f)(t)=\sum_{m=0}^{\infty}f(t-m\log p).
\end{equation}
These operators are related by an \emph{Euler product} as follows:

\begin{equation}
f(t)\mapsto \mathfrak{a}(f)(t)=\left(\prod_{p\in \mathcal{P}}\mathfrak{a}_{p}\right)(f)(t),\label{Eq:EPr}
\end{equation} 
where the product is the composition of operators.\\

\hspace*{3mm}Let $f$ be an infinitely differentiable function on $\mathbb{R}$.\,Then, the Taylor series of $f$ can be formally written as
\begin{align}
f(t+h)&=f(t)+\frac{hf'(t)}{1!}+\frac{h^{2}f^{''}(t)}{2!}+...\notag\\
      &=e^{h\frac{d}{dt}}(f)(t)=e^{h\partial}(f)(t),\label{Eq:Taylexp}
\end{align}       
where $\partial=\frac{d}{dt}$ is the first order differential operator with respect to $t$.\footnote{This differential operator is the infinitesimal generator of the (one-parameter) group of shifts on the real line.\,For this reason, it is also called the \emph{infinitesimal shift}; see Lemma \ref{Lem:contsem} and Lemma \ref{Lem:Infsht} which justify this terminology.}

\begin{remark}\label{Rk:Ctfn}
In our later, more mathematical discussion, $f$ will not necessarily be the counting function of some generalized fractal string $\eta$, but will instead be allowed to be an element of the Hilbert space $\mathbb{H}_{c}$ $($with possibly some additional conditions on $f$ or on the parameter $c$$)$;  see Equation \emph{(\ref{Eq:HS})} below and the text surrounding it, along with Equations \emph{(\ref{Eq:acf})} and \emph{(\ref{Eq:Dsp})}.
\end{remark}

\hspace*{3mm}Note that this yields a new \emph{heuristic} representation for the spectral operator and its prime factors:

\begin{align}
\mathfrak{a}(f)(t)&=\sum_{n=1}^{\infty}e^{-(\log n)\partial}(f)(t)=\sum_{n=1}^{\infty}\left(\frac{1}{n^{\partial}}\right)(f)(t)\notag\\
                  &=\zeta(\partial)(f)(t)=\zeta_{\mathfrak{h}}(\partial)(f)(t)=\prod_{p\in \mathcal{P}}(1-p^{-\partial})^{-1}(f)(t)\label{Eq:Spr}
\end{align}  
and for any prime $p$,
\begin{align}
\mathfrak{a}_{p}(f)(t)&=\sum_{m=0}^{\infty}f(t-m\log p)=\sum_{m=0}^{\infty}e^{-m(\log p)\partial}(f)(t)=\sum_{m=0}^{\infty}\left(p^{-\partial}\right)^{m}(f)(t)\notag\\
                      &=\left(\frac{1}{1-p^{-\partial}}\right)(f)(t)=(1-p^{-\partial})^{-1}(f)(t)=\zeta_{\mathfrak{h}_{p}}(\partial)(t).\label{Eq:Pf}
\end{align}

\begin{remark}The above representations of the spectral operator, its operator-valued Euler factors and its operator-valued Euler product were given in \emph{[\textbf{La-vF3}, \S6.3.2]} without specifying a domain $($or a `core'$)$ enabling one to study them and analyze some of their fundamental properties.\,$($See footnote $12$$.)$\,Finding an appropriate Hilbert space and a domain which is equipped with natural boundary conditions satisfied by the class of counting functions of generalized fractal strings was one of the first steps taken in \emph{\cite{HerLa1}} prior to studying these operators and then deriving the desired spectral reformulation of the Riemann hypothesis.\,Finally, we mention the fact that the operator-valued prime factors and their operator-valued Euler product are investigated in \emph{\cite{HerLa2}}.\,In particular, in that paper, we establish the convergence $($in the operator norm$)$ of the operator-valued Euler product, when $c>1$, and investigate the conjecture $($suggested by comments in $[$\emph{\textbf{La-vF3}}, \S6.3.2$]$$)$ according to which, in an appropriate sense, this same Euler product can be analytically continued and shown to converge to the spectral operator even in the critical strip $0<Re(s)<1$ $($i.e., for $0<c<1$$)$.  
\end{remark}

\subsection{The weighted Hilbert space $\mathbb{H}_{c}$.}
\hspace*{3mm}In \cite{HerLa1},  we start by providing a functional analytic framework enabling us to rigorously study the spectral operator.\,This functional analytic framework is based in part on defining a specific weighted Hilbert space $\mathbb{H}_{c}$, depending on a parameter $c\geq 0$, in which the spectral operator is acting, and then on precisely defining and studying this operator.\,We set
\begin{equation}\label{Eq:HS}
\mathbb{H}_{c}=L^{2}(\mathbb{R},\mu_{c}(dt)),
\end{equation}
\text
where $\mu_{c}$ is the absolutely continuous measure on $\mathbb{R}$ given by $\mu_{c}(dt):=e^{-2ct}dt$ (here, $dt$ is the Lebesgue measure on $\mathbb{R}$).\\

\begin{remark}
Note that $\mathbb{H}_{c}$ is the space of $($$\mathbb{C}$-valued$)$ Lebesgue square-integrable functions $f$ with respect to the positive \emph{weight function} $w(t)=e^{-2ct}$$:$
\begin{equation}\label{Eq:fc}
||f||^{2}_{c}:=\int_{\mathbb{R}}|f(t)|^{2}e^{-2ct}dt<\infty.
\end{equation}
\text
It is obtained by completing the space $\mathcal{H}_{c}$ of infinitely differentiable functions $f$ on $\mathbb{R}=(-\infty,+\infty)$ satisfying the finiteness condition $($\ref{Eq:fc}$)$.\,$($It follows, of course, that $\mathcal{H}_{c}$ is dense in $\mathbb{H}_{c}$.$)$
\end{remark}

The Hilbert space $\mathbb{H}_{c}$ is equipped with the inner product
\[<f,g>_{c}\,:=\int_{\mathbb{R}}f(t)\overline{g(t)}e^{-2ct}dt\]
and the associated Hilbert norm $||.||_{c}=\sqrt{<.\,,.>_{c}}$ (so that $||f||_{c}^{2}=$\\$\int_{\mathbb{R}}|f(t)|^{2}e^{-2ct}dt)$.\,Here, $\overline{g}$ denotes the complex conjugate of $g$.\\

\hspace*{3mm}Next, we introduce the boundary conditions which are naturally satisfied within our framework by the class of counting functions of generalized fractal srings.\,(See Remarks \ref{Rk:Ctfn}, \ref{Rk:ascond} and \ref{Rk:2.7}.)\,Note that if $f\in \mathbb{H}_{c}$ and $f$ is absolutely continuous on $\mathbb{R}$ (i.e., $f\in AC(\mathbb{R})$), then 
\begin{equation}
|f(t)|e^{-ct}\to 0 \mbox{\quad as\,\,$t\rightarrow \pm \infty$}, \label{Eq:Bdcond}
\end{equation}
respectively.\,Because the domain $D(\partial_{c})$ of the infinitesimal shift $\partial_{c}$ will consist of absolutely continuous functions (see Equation (\ref{Eq:acf})), these are \emph{natural boundary conditions}, in the sense that they are automatically satisfied by any function $f$ in the domain of $\partial_{c}$ or of a function of $\partial_{c}$, such as the spectral operator $\mathfrak{a}_{c}=\zeta(\partial_{c})$ (see Equation (\ref{Def:Spop})).\,Furthermore, observe that the boundary conditions (\ref{Eq:Bdcond}) imply that $|f(t)|=o(e^{-c|t|})$ as $t\to -\infty$ and hence (since $c\geq 0$), that $f(t)\to 0$ as $t\to -\infty$.\\

\begin{remark}\label{Rk:ascond}
The asymptotic condition at $+\infty$ in Equation \emph{(\ref{Eq:Bdcond})} implies that $($roughly speaking$)$ the functions $f$ satisfying these \emph{boundary conditions} correspond to elements of the space of fractal strings with dimension $D\leq c$; see also Remark \ref{Rk:2.7}.
\end{remark}

\begin{remark}\label{Rk:2.7}
Note that in the original multiplicative variable $x=e^{t}$ and for an ordinary fractal string $\mathcal{L}$, it is shown in \emph{\cite{LaPo2}} that if $f(t):=N_{\mathcal{L}}(x)$ is of order not exceeding $($respectively, is precisely of the order of$)$ $x^{c}=e^{ct}$ as $x\to +\infty$ $($i.e., as $t\to +\infty$$)$, then $D\leq c$ $($respectively, $c=D$, the Minkowski dimension of $\mathcal{L}$$)$.\footnote{See \cite{LaPo2} (and \S4.2 below) for a thorough discussion of the geometric and spectral interpretations of various asymptotic conditions satisfied by the counting functions of ordinary fractal strings.\,(See also \cite{HeLa} for further generalizations.)}\,In addition, it follows from \emph{\cite{La2, LaPo2, La-vF3, LaLu-vF1}} that $($with the same notation as above$)$
\begin{equation}
D=\alpha:=\inf\{\gamma \geq 0:\,N_\mathcal{L}(x)=f(e^{t})=O(e^{\gamma t}),\,\mbox{\emph{as}}\,\,t\to +\infty\},
\end{equation}
\text
and hence, $\alpha$ coincides with the abscissa of convergence $D=D_{\mathcal{L}}$ of the geometric zeta function $\zeta_{\mathcal{L}}$.\\

\hspace*{3mm}Moreover, let us suppose that $\mathcal{L}$ is normalized so that its geometric counting function satisfies $N_{\mathcal{L}}(x)=0$ for $0<x\leq 1$ $($which, in the additive variable $t=\log x$, amounts to assuming that $f(t)= 0$ for all $t\leq 0$, where we let $f(t):=N_{\mathcal{L}}(e^{t})$, as above$)$.\footnote{Without loss of generality, this can always be done since there exists $x_{0}>0$ such that $N_{\mathcal{L}}(x)=0$ for all $0<x\leq x_{0}$.\,Indeed, it suffices to replace each $l_{j}$ with $\frac{l_{j}}{l_{1}}$ to allow the choice $x_{0}=1$ (in the multiplicative variable, and hence, $t=0$, in the additive variable).}\,Then we can simply let $F(t):=f(t)$ for $t\geq 0$ and $F(t):=0$ for $t\leq 0$ in order to obtain a nonnegative function $F$ defined on all of $\mathbb{R}$, vanishing identically on $(-\infty,0]$, and having the same asymptotic behavior as $f(t)$ as $t\to +\infty$.\,In particular, if $f\in L^{2}([0,+\infty), \mu_{c}(dt))$ satisfies $f(t)=o(e^{ct})$ as $t \to +\infty$, then $F\in \mathbb{H}_{c}$ and satisfies the above boundary conditions stated in Equation $($\ref{Eq:Bdcond}$)$$:\,F(t)=o(e^{ct})$ as $t \to \pm \infty$.\,Note that if, furthermore, $f$ is absolutely continuous on $[0,+\infty)$ $($i.e., $f\in AC([0,+\infty))$$)$, then $F$ is absolutely continuous on all of $\mathbb{R}$ $($because its almost everywhere defined derivative $F'$ is locally integrable on $\mathbb{R}$ and $F(t)=\int_{0}^{t}F'(\tau)d\tau$ for all $t\in\mathbb{R}$$)$ and hence belongs to the domain of $\partial_{c}$, as defined by Equation \emph{(\ref{Eq:acf})} below.\,In  particular, as was noted earlier, $F$ then automatically satisfies the boundary condition at $+\infty$.
\end{remark}
\subsection{The infinitesimal shifts $\partial_{c}$ and their properties.}

\hspace*{3mm}In this subsection, we first define the domain of the infinitesimal shift $\partial_{c}$, in \S3.1.1, then review the properties of $\partial_{c}$  (and of its spectrum) established in \cite{HerLa1}, in \S3.3.2, and finally study (in \S3.3.3) the contraction group of linear operators generated by $\partial_{c}$; as it turns out, this is a suitable version of the shift group on the real line.\\

\subsubsection{The domain of the infinitesimal shifts}
\hspace*{3mm}Recall from the heuristic discussion surrounding Equation (\ref{Eq:Spr}) that the differential operator $\partial=\partial_{c}$, also called the \emph{infinitesimal shift}, arises naturally in the representation of the spectral operator, its operator-valued Euler factors and its operator-valued Euler product (see Equation (\ref{Eq:EPr})).\,Motivated by this fact, and in light of our proposed definition for the spectral operator in Equation (\ref{Def:Spop}), we adopt the following precise domain for the infinitesimal shift $\partial_{c}$:
\begin{equation}
D(\partial_{c})=\{f\in \mathbb{H}_{c}\cap AC(\mathbb{R}):\,f'\in \mathbb{H}_{c}\},\label{Eq:acf}
\end{equation}
\text
where $AC(\mathbb{R})$  is the space of (locally) absolutely continuous functions on $\mathbb{R}$ and $f'$ denotes the derivative of $f$, viewed either as a function or a distribution.\,Recall that for $f\in AC(\mathbb{R})$, $f'$ exists pointwise almost everywhere and is locally integrable on $\mathbb{R}$, therefore defining a regular distribution on $\mathbb{R}$.\,(See, e.g., \cite{Br, Fo, Ru, Schw}.)\footnote{Note that $D(\partial_{c})$ can be viewed as the weighted Sobolev space $H^{1}(\mathbb{R},\mu_{c}(dt))$.\,See Remark \ref{partSob}; see also, e.g., \cite{Br} or \cite{Fo} for the classic case when $c=0$ and hence this space coincides with the standard Sobolev space $H^{1}(\mathbb{R})$.}\,In addition, we let
\begin{equation}
\partial_{c}:=f',\mbox{\quad for\,} f\in D(\partial_{c}).
\end{equation}

\begin{remark}\label{partSob}
Alternatively, one could view $\partial_{c}$ as a bounded $($normal$)$ linear operator acting on the weighted Sobolev space $D(\partial_{c})$, equipped with the Hilbert norm $N_{c}(f):=(||f||_{c}^{2}+||f'||_{c}^{2})^{\frac{1}{2}}$.\,We will not abopt this point of view here, although it is helpful in order to motivate some of the proofs of the results obtained in \cite{HerLa1}.
\end{remark}

\subsubsection{Normality and spectra of the infinitesimal shifts.}
\hspace*{3mm}Our first result will enable us to form various functions of the first order differential operator $\partial_{c}$ and, in particular, to precisely define the spectral operator $\mathfrak{a}_{c}$.

\begin{theorem}\label{Thm:part}
\emph{\cite{HerLa1}}\,$\partial_{c}$ is an unbounded normal\footnote{Recall that this means that $\partial_{c}$ is a closed $($and densely defined$)$ operator which commutes with its adjoint $\partial^{*}_{c}$; see \cite{Ru}.} linear operator on $\mathbb{H}_{c}$.\\Moreover, its adjoint $\partial_{c}^{*}$ is given by 
\begin{equation}\label{Eq:part}
\partial_{c}^{*}=2c-\partial_{c},\mbox{\quad\emph{with domain}\,}\,D(\partial_{c}^{*})=D(\partial_{c}).
\end{equation}
\end{theorem}

\begin{remark}
We encourage the reader to consult  Appendix B for a sketch of the proof of  Theorem \ref{Thm:part} and for a useful reformulation of that theorem, provided in Corollary \ref{CorTh}.\,This proof and the corresponding reformulation are  based on a representation of the infinitesimal shift $\partial_{c}$ in terms of a linear unbounded normal operator $V_{c}$ $($acting on $\mathbb{H}_{c}=L^{2}(\mathbb{R}, e^{-2ct}dt)$$)$ which is unitarily equivalent to the standard momemtum operator $V_{0}=\frac{1}{i}\partial=\frac{1}{i}\frac{d}{dt}$ $($acting on $\mathbb{H}_{0}=L^{2}(\mathbb{R})$$)$.\,Namely, we have $\partial_{c}=c+i V_{c}$.
\end{remark}

\hspace*{3mm}In order to find the spectrum $\sigma(\mathfrak{a}_{c})$ of the spectral operator, we first determine the spectrum of $\partial_{c}$, which turns out to be equal to the vertical line of the complex plane passing through the constant $c$.

\begin{theorem}\label{Thm:spp}
\emph{\cite{HerLa1}}\,Let $c\geq 0$.\,Then, the spectrum of $\partial_{c}$ is the closed vertical line of the complex plane passing through $c$.\,Furthermore, it coincides with the essential spectrum, $\sigma_{e}(\partial_{c})$, of $\partial_{c}:$
\begin{equation}\label{Eq:spc}
\sigma(\partial_{c})=\sigma_{e}(\partial_{c})=\{\,\lambda\in \mathbb{C}:\,Re(\lambda)=c\,\}.
\end{equation}

\hspace*{3mm}More specifically, the point spectrum of $\partial_{c}$ is empty $($i.e., $\partial_{c}$ does not have any eigenvalues$)$\footnote{We caution the reader that (perhaps surprisingly) the terminology concerning the spectra of unbounded operators is not uniform throughout the well-developed literature on this classical subject; see, e.g., \cite{DunSch, EnNa, Kat, ReSi, Sc, JoLa}.} and $\sigma_{ap}(\partial_{c})$, the \emph{approximate point spectrum} of $\partial_{c}$, coincides with $\sigma(\partial_{c})$.\,Hence, $\sigma_{ap}(\partial_{c})$ is also given by the right-hand side of Equation \emph{(\ref{Eq:spc})}.\footnote{Recall that $\lambda \in \sigma_{ap}(\partial_{c})$ (i.e., $\lambda$ is an \emph{approximate eigenvalue} of $\partial_{c}$) if and only if there exists a sequence $\{f_{n}\}^{\infty}_{n=1}$ of elements of $D(\partial_{c})$ such that $||f_{n}||=1$ for all $n\geq 1$ and $||\partial_{c}f_{n}-\lambda f_{n}||_{c}\to 0$ as $n\to \infty$. (See, e.g., \cite{EnNa} or \cite{Sc}.)}
\end{theorem}

\begin{remark}
Note that the infinitesimal shift $\partial_{c}$ is unbounded $($since, by Theorem \ref{Thm:spp}, its spectrum is unbounded$)$, normal $($by Theorem \ref{Thm:part}, $\partial^{*}_{c}\partial_{c}=\partial_{c}\partial^{*}_{c}$$)$, and sectorial $($in the extended sense of \emph{\cite{Ha}}, since by Theorem \ref{Thm:spp}, $\sigma(\partial_{c})$ is contained in a sector of angle $\frac{\pi}{2}$$)$.
\end{remark}

\begin{remark}
It is shown in \emph{\cite{HerLa1}} that in addition to being normal, $\partial_{c}$ is $m$-accretive, in the sense of \emph{\cite{Kat}}.\,According to a well-known theorem in semigroup theory, this means that $\partial_{c}$ is the infinitesimal generator of a contraction semigroup of operators; see Lemmas \ref{Lem:contsem} and \ref{Lem:Infsht}.
\end{remark}

\subsubsection{The strongly continuous group of operators $\{e^{-t\partial_{c}}\}_{t\in \mathbb{R}}$}
\hspace*{3mm}The strongly continuous contraction group\footnote{We refer to \cite{EnNa, Go, HiPh, JoLa, Kat, Paz, ReSi} for the theory of strongly continuous semigroups.} of bounded linear operators $\{e^{-t\partial_{c}}\}_{t\in \mathbb{R}}$ plays a crucial role in the representation of the spectral operator $\mathfrak{a}_{c}=\zeta(\partial_{c})$ which was obtained and rigorously justified in \cite{HerLa1}; see Theorem \ref{Thm:11} and Remark \ref{RK:sden}.\,(See also Equations (\ref{Eq:Taylexp}) and (\ref{Eq:Spr}), along with the discussion surrounding them, concerning heuristic representations of the spectral operator, its operator-valued Euler factors and its operator-valued Euler product.)\\

\hspace*{3mm} Using Theorem \ref{Thm:part} (or equivalently, Corollary \ref{CorTh} of Appendix B), we obtain the following result.

\begin{lemma}\emph{[\textbf{HerLa1}]}\label{Lem:contsem}
For any $c\geq 0$, $\{e^{-t\partial_{c}}\}_{t\in \mathbb{R}}$ is a strongly continuous contraction group of operators and $||e^{-t\partial_{c}}||= e^{-tc}$ for any $t\in\mathbb{R}$.\,The adjoint group $\{(e^{-t\partial_{c}})^{*}\}_{t\in \mathbb{R}}$ is then given by $\{e^{-t\partial_{c}^{*}}\}_{t\in\mathbb{R}}=\{e^{-t(2c-\partial_{c})}\}_{t\in\mathbb{R}}$.
\end{lemma}

\begin{remark}\label{Rk:utcor}
If we let $i:=\sqrt{-1}$, it follows from Corollary \ref{CorTh} that $\{e^{-t(\frac{\partial_{c}}{i})}\}_{t\in\mathbb{R}}$ is a \emph{unitary group} if and only if $c=0$.\,$($Compare with Theorem \ref{Thm:mom} in Appendix B.$)$
\end{remark}

Another key feature of this strongly continous group of operators is highlighted in the following result.
\begin{lemma}\emph{[\textbf{HerLa1}]}\label{Lem:Infsht} 
For any $c\geq 0$, the strongly continuous group of operators $\{e^{-t\partial_{c}}\}_{t\in \mathbb{R}}$ is a translation $($or \emph{shift}$)$ group.\,That is, for every $t\in \mathbb{R}$, $(e^{-t\partial_{c}})(f)(u)=f(u-t)$, for all $f\in \mathbb{H}_{c}$ and $u\in \mathbb{R}$.\,$($For a fixed $t\in \mathbb{R}$, this equality holds between elements of $\mathbb{H}_{c}$ and hence, for a.e. $u\in\mathbb{R}$$.)$
\end{lemma}

In light of Lemma \ref{Lem:Infsht}, the infinitesimal generator $\partial=\partial_{c}$ of the shift group $\{e^{-t\partial}\}_{t\in \mathbb{R}}$ is called the \emph{infinitesimal shift} of the real line.\\

\subsection{The spectral operator $\mathfrak{a}_{c}$.}
\hspace*{3mm}In \cite{HerLa1}, we define the \textit{spectral operator} $\mathfrak{a}=\mathfrak{a}_{c}$ as follows, where $\partial=\partial_{c}$ is the normal operator defined in \S3.3.1:
\begin{equation}\label{Def:Spop}
\mathfrak{a}=\zeta(\partial),
\end{equation}
via the measurable functional calculus for unbounded normal operators; see, e.g., \cite{Ru}.\,If, for simplicity, we assume $c\ne 1$ to avoid the pole of $\zeta$ at $s=1$, then (in light of Theorem \ref{Thm:spp}) $\zeta$ is holomorphic (and, in particular, continuous) in an open neighborhood of $\sigma(\partial)$.\,If $c=1$ is allowed, then (still by Theorem \ref{Thm:spp}) $\zeta$ is meromorphic in an open neighborhood of $\sigma(\partial)$ (actually, in all of $\mathbb{C}$).\,Hence, when $c\ne 1$, we could simply use the holomorphic (or the continuous) functional calculus for unbounded normal operators (see \cite{Ru}), whereas when $c=1$, we could use the meromorphic functional calculus for sectorial operators\,(see \cite{Ha}).\,For any value of $c$, however, the measurable functional calculus can be used.\\

\hspace*{3mm}The domain of the spectral operator is the following:
\begin{equation}\label{Eq:Dsp}
D(\mathfrak{a})=\{f\in D(\partial):\,\mathfrak{a}(f)=\zeta(\partial)(f)\in \mathbb{H}_{c}\}.
\end{equation}

\hspace*{3mm}Our next result, Theorem \ref{Thm:11} below, provides a representation of the spectral operator $\mathfrak{a}$ as a composition map of the Riemann zeta function $\zeta$ and the first order differential operator $\partial_{c}$.\,It also gives a natural connection between this representation and the earlier one obtained for the spectral operator in Equations (\ref{Eq:spop}) and (\ref{Eq:Spr}).\,(See also Lemmas \ref{Lem:contsem}, \ref{Lem:Infsht} and  Equations (\ref{Eq:EPrd}), (\ref{Eq:3.10}) in \S7.)

\begin{theorem}\emph{\cite{HerLa1}}\,Assume that $c>1$.\,Then, $\mathfrak{a}$ can be uniquely extended to a bounded operator on $\mathbb{H}_{c}$ and, for any $f\in \mathbb{H}_{c}$, we have $($for almost all $t\in \mathbb{R}$ or as an equality in $\mathbb{H}_{c}$$):$\label{Thm:11}
\begin{equation}
\mathfrak{a}(f)(t)=\sum_{n=1}^{\infty}f(t-\log n)=\zeta(\partial_{c})(f)(t)=\left(\sum_{n=1}^{\infty}n^{-\partial_{c}}\right)(f)(t).\label{Eq:spcf}
\end{equation}
\end{theorem}
In other words, for $c>1$, we have
\begin{equation}
\mathfrak{a}_{c}=\zeta(\partial_{c})=\sum_{n=1}^{\infty}n^{-\partial_{c}},\label{Eq:2.23}
\end{equation}
\text
where the equality holds in $\mathcal{B}(\mathbb{H}_{c})$, the space of bounded linear operators on $\mathbb{H}_{c}$.\\

\begin{remark}\label{RK:sden}
For any $c>0$, we also show in \cite{HerLa1} that Equation $($\ref{Eq:spcf}$)$ holds for all $f$ in a suitable dense subspace of $D(\mathfrak{a})$, which we conjectured to be a core for $\mathfrak{a}$ and hence to uniquely determine the unbounded operator $\mathfrak{a}=\mathfrak{a}_{c}=\zeta(\partial_{c})$, viewed as the $($operator-valued$)$ `\emph{analytic continuation}' of \,$\sum_{n=1}^{\infty}n^{-\partial_{c}}$ to the critical strip $0<Re(s)<1$ $($and thus also to the open half-plane $Re(s)>0$$)$.
\end{remark}

\hspace*{3mm}In order to study the invertibility of the spectral operator, a characterization of the spectrum $\sigma(\mathfrak{a}_{c})$ of the spectral operator was obtained in \cite{HerLa1} by using the spectral mapping theorem for unbounded normal operators (the continuous version when $c\ne 1$ and the meromorphic version, when $c=1$); see Remark \ref{Rk:spt}.

\begin{theorem}\emph{\cite{HerLa1}}\,Assume that $c\geq0$.\,Then \label{Thm:ssop}
\begin{equation}\label{Eq:Merv}
\sigma(\mathfrak{a})=\overline{\zeta(\sigma(\partial))}=cl\big(\zeta(\{\lambda\in\mathbb{C}:\,Re(\lambda)=c\})\big),
\end{equation}
where $\sigma(\mathfrak{a})$\,is the spectrum of $\mathfrak{a}=\mathfrak{a}_{c}$ and $\overline{N}=cl(N)$ is the closure of $N\subset \mathbb{C}$.
\end{theorem}

\begin{remark}\label{Rk:spt}
We refer to the appropriate appendix in \cite{HerLa1} for a discussion of the spectral mapping theorem for linear unbounded normal operators.\,In short, if $\phi$ is a continuous function on $\sigma(\mathcal{Q})$, where $\mathcal{Q}$ is a given $($possibly unbounded$)$ normal operator, then $\sigma(\phi(\mathcal{Q}))=\overline{\phi(\sigma(\mathcal{Q}))}$.\,Moreover, if $\phi$ is a $($$\mathbb{C}$-valued$)$ meromorphic function on an open neighborhood of the spectrum $\sigma(\mathcal{Q})$ $($and say, $\mathcal{Q}$ has no eigenvalues$)$,\footnote{Or more generally, if no eigenvalue of $\mathcal{Q}$ coincides with a pole of $\phi$ lying in $\sigma(\mathcal{Q})$.} then $\tilde{\sigma}(\phi(\mathcal{Q}))=\phi(\tilde{\sigma}(\mathcal{Q}))$, interpreted as an equality between subsets of the Riemann sphere $\tilde{\mathbb{C}}:=\mathbb{C}\cup\{\infty\}$.\,Here, given a linear operator $K$, the extended spectrum $\tilde{\sigma}(K)$ is the compact subset of $\tilde{\mathbb{C}}$ defined by $\tilde{\sigma}(K)=\sigma(K)$ if $K$ is bounded and $\tilde{\sigma}(K)=\sigma(K)\cup\{\infty\}$ if $K$ is unbounded.\,$($Note that if $\phi$ is meromorphic, then it is continuous when viewed as a $\tilde{\mathbb{C}}$-valued function.$)$\footnote{We wish to thank Daniel Lenz and Markus Haase for helpful written correspendences about the spectral mapping theorem in this context.}
\end{remark}

\hspace*{3mm}We will see in \S5 that the characterization of the spectrum of the infinitesimal shift $\partial_{c}$ obtained in Theorem \ref{Thm:spp} will play an important role in our proposed definition of the \emph{truncated infinitesimal shifts} and \emph{spectral operators}, $\{\partial^{(T)}_{c}\}_{T\geq0}$ and $\{\mathfrak{a}^{(T)}_{c}\}_{T\geq0}$ (in \S5.1), the determination of their corresponding spectra $\sigma(\partial_{c}^{(T)})$ and $\sigma\big(\mathfrak{a}^{(T)}_{c}\big)$ (in \S5.2), the study of the quasi-invertibility of $\mathfrak{a}_{c}$ (in \S5.3), and ultimately, in our spectral reformulation of the Riemann hypothesis discussed in \S6.1.

\begin{center}
\section{Inverse and Direct Spectral Problems for Fractal Strings}
\end{center}
\subsection{The original inverse spectral problem.}
\hspace*{3mm}The problem of deducing geometric information from the spectrum of a fractal string, or equivalently, of addressing the question
\begin{quotation}
\hspace*{1.5 cm}\emph{\textquotedblleft Can one hear the shape of a fractal string?\textquotedblright},
\end{quotation}
was first studied by the second author and H.\,Maier in \cite{LaMa1,LaMa2}.\,More specifically, \emph{the inverse spectral problem} they considered was the following: \\

\textquotedblleft \emph{Given any fixed} $D\in (0,1)$, \emph{and any fractal string} $\mathcal{L}$ \emph{of dimension D such that for some constant} $c_{D}>0$ \emph{and} $\delta>0$, \emph{we have}
\begin{equation}\label{Eq:Wt}
N_{\nu}(x)=W(x)-c_{D}x^{D}+ O(x^{D-\delta}),\mbox{\quad as\,\,} x\to +\infty,
\end{equation}
\emph{is it true that} $\mathcal{L}$ \emph{is Minkowski measurable}?\textquotedblright.\\
\begin{remark}\label{Rk:W}
Here, the Weyl term $W(x)$ is the leading asymptotic term.\,Namely,
\begin{equation}
W(x):=(2\pi)^{-1}vol_{1}(\Omega)x,\label{Eq:Weylter}
\end{equation}
where $x$ is the $($normalized$)$ frequency variable and $vol_{1}(\Omega)$ denotes the \textquotedblleft volume\textquotedblright\\$($really, the length$)$ of $\Omega\subset \mathbb{R}$.\,Furthermore, much as before, the spectral counting function $N=N_{\nu}(x)$ is equal to the number of frequencies of $\mathcal{L}$ less than $x$.
\end{remark}
\hspace*{3mm}The geometric notion of Minkowski measurability will be recalled below in Definition \ref{Def:Minmeas}.\,For now, we note that the above question is indeed stated in the form of an inverse spectral problem.\,Namely, one is asked to deduce geometric information about a fractal string from spectral asymptotic information about the string.\,Roughly speaking, given that the spectrum of $\mathcal{L}$ has a monotonic asymptotic second term (i.e., does not have any oscillations of order $D$, the Minkowski dimension of $\mathcal{L}$ (or $\partial\Omega$) (see Equation (\ref{Eq:Wt}) and Definition \ref{Def:Minmeas}), does it follow that the geometry of $\mathcal{L}$ does not have any oscillations of leading order $D$ (see Equation (\ref{Eq:Mm}) in Definition \ref{Def:Minmeas}, along with \S4.2 below)?\\

\hspace*{3mm}The authors of \cite{LaMa1, LaMa2} have shown that this question \emph{\`a la Mark Kac} (but interpreted rather differently than in [\textbf{Kac}]) \textquotedblleft Can one hear the shape of a fractal string?\textquotedblright, is intimitely connected with the Riemann hypothesis.\,More specifically, they proved that for a given $D\in(0,1)$, this inverse spectral problem is true for every fractal string of dimension $D$ if and only if the Riemann zeta function does not have any zeroes along the vertical line $Re(s)=D$:\,$\zeta(s)\ne 0$ for $Re(s)=D$.\\

\hspace*{3mm}It follows, in particular, that the inverse spectral problem has a negative answer in the `mid-fractal case' where $D=\frac{1}{2}$ (because $\zeta$ has a zero, and even infinitely many zeroes, along the critical line $Re(s)=\frac{1}{2}$).\,Moreover, it follows that this inverse spectral problem has a positive answer for all fractal strings whose dimension is an arbitrary number $D\in(0,1)-\frac{1}{2}$ if and only if the Riemann hypothesis is true.
\begin{remark}
The work in \emph{\cite{LaMa2}} $($announced in \cite{LaMa1}$)$ was revisited and extended to a large class of arithmetic zeta functions in \emph{\cite{La-vF1, La-vF2, La-vF3}}, using the explicit formulas recalled in Theorem \ref{Thm:2} and the then rigorously defined notion of complex dimension.\,$($See \emph{[\textbf{La-vF3}, Ch.\,9]}.$)$ Our work in \emph{\cite{HerLa1, HerLa2}} can also be extended to this more general setting, as will be clear to the reader familiar both with the functional calculus $($for unbounded normal operators$)$ and the theory of $L$-functions, but by necessity of concision, we will not discuss this development here.
\end{remark}

\begin{definition}\label{Def:Minmeas}
A fractal string $\mathcal{L}$ $($or equivalently, $\partial\Omega$, the boundary of the associated open set $\Omega\subset \mathbb{R}$$)$ is said to be \emph{Minkowski measurable} if the following limit exists in $(0,+\infty):$
\begin{equation}
\lim_{\epsilon\to 0^{+}}\frac{vol_{1}(\Omega_{\epsilon})}{\epsilon^{1-D}}:=\mathcal{M}(\mathcal{L}),\label{Eq:Mm}
\end{equation}
where $\mathcal{M}(\mathcal{L})$ is called the \emph{Minkowski content} of $\mathcal{L}$ $($or of $\partial{\Omega}$$)$ and, as before, $vol_{1}(\Omega_{\epsilon})$ denotes the volume of the inner $\epsilon$-neighborhood of $\partial \Omega:$ $\Omega_{\epsilon}=\{x\in \Omega:\,dis(x,\partial \Omega)<\epsilon\}$.\,It then follows that $D$ is the Minkowski $($or box$)$ dimension of $\mathcal{L}$ $($i.e., of $\partial \Omega$$)$.\\
\hspace*{3mm}Moreover, recall that the \emph{Minkowski dimension} of $\mathcal{L}$ $($or equivalently, of $\partial\Omega$$)$ is defined by
\begin{equation}
D=\sup\{\alpha\geq0:\mathcal{M}^{*}_{\alpha}(\mathcal{L})=\infty\}=\inf\{\alpha\geq0:\mathcal{M}^{*}_{\alpha}(\mathcal{L})=0\},
\end{equation}
where $\mathcal{M}^{*}_{\alpha}(\mathcal{L})$, the $\alpha$-\emph{dimensional upper Minkowski content} of $\mathcal{L}$ $($or of $\partial\Omega$$)$, is given by
\begin{equation}\label{Eq:upMink}
\mathcal{M}^{*}_{\alpha}(\mathcal{L}):=\limsup_{\epsilon\to 0^{+}}\frac{vol_{1}(\Omega_{\epsilon})}{\epsilon^{1-\alpha}}.
\end{equation}
$($The $\alpha$-\emph{dimensional lower Minkowski content} of $\mathcal{L}$, $\mathcal{M}_{*,\,\alpha}(\mathcal{L})$, is defined analogously, but with a lower limit instead of an upper limit on the right-hand side of the counterpart of Equation $($\ref{Eq:upMink}$)$.$)$
\end{definition}
\begin{remark}\label{RK:MDabcv}
Recall that $D=D_{\mathcal{L}}$, the Minkowski dimension of a fractal string $\mathcal{L}$, coincides with the abscissa of convergence of the Dirichlet series initially defining the geometric zeta function $\zeta_{\mathcal{L}}$; see Equation $($\ref{Eq:Dimeta}$)$ and the text following it.\footnote{This is why we abuse notation by using the same symbol for these two notions.}\,This key fact was first observed in \cite{La2} using an important result of Besicovitch and Taylor \emph{\cite{BesTa}}, and a direct proof of this equality was later provided in \emph{[\textbf{La-vF2}, Thm.\,1.10]}.\,Furthermore, recall that in the present geometric situation, we always have $D\in[0,1]$.\, In other words, the dimension $D=D_{\mathcal{L}}$ of an ordinary fractal string always lies in the `\emph{critical interval}' $(0,1)$ or coincides with one of its endpoints, $0$ and $1$, corresponding to the `least' and `most' fractal case, respectively $($in the terminology of \cite{La1}$)$.
\end{remark}

\subsection{Fractal strings and the (modified) Weyl--Berry conjecture.}

\hspace*{3mm}Prior to the work in \cite{LaMa1, LaMa2}, the second author and Carl Pomerance \cite{LaPo1, LaPo2} had studied the corresponding \emph{direct spectral problem} for fractal strings.\\They thereby had resolved in the affirmative the (one-dimensional) modified Weyl--Berry conjecture (as formulated in \cite{La1}) according to which if a fractal string $\mathcal{L}$ is Minkowski measurable of dimension $D\in(0,1)$, then its spectral counting function $N_{\nu}(x)$ has a monotonic asymptotic second term, proportional to $\mathcal{M}(\mathcal{L})x^{D}$.\,More specifically, the authors of \cite{LaPo1, LaPo2} had shown, in particular, that if $\mathcal{L}$ is Minkowski measurable $($which, according to a key result in \cite{LaPo2}, is true iff $l_{j}\sim L.j^{-\frac{1}{D}}$ as $j\to \infty$ or equivalently, iff $N_{\mathcal{L}}(x)\sim C.x^{D}$ as $x\to +\infty$, for some $L>0$ and $C>0$$)$,\footnote{Here, $l_{j}\sim m_{j}$ as $j\to \infty$ means that $l_{j}=m_{j}(1+o(1))$ as $j\to\infty$, where $o(1)$ stands for a function tending to zero at infinity; and similarly for functions of a continuous variable $x\in(0,+\infty)$.} then the eigenvalue $($or rather, frequency$)$ counting function $N_{\nu}(x)$ satisfies Equation (\ref{Eq:Wt}) $($with $o(x^{D})$ instead of $O(x^{D-\delta})$$)$, where
\begin{equation}
c_{D}:=2^{-(1-D)}(1-D)(-\zeta(D))\mathcal{M}(\mathcal{L})
\end{equation}
\text
and $\mathcal{M}(\mathcal{L})$ is the Minkowski content of $\mathcal{L}$, as defined in Equation (\ref{Eq:Mm}).\,$($Note that $-\zeta(D)>0$ for $0<D<1$.$)$\\

\hspace*{3mm}This result and its proof (along with related results and conjectures in \cite{La1} and \cite{La2, La3}) suggested the possibility of developing a theory of complex exponents (or complex dimensions) which enables one to give a natural geometric meaning to the critical strip $0<Re(s)<1$.\,Accordingly, the least (respectively, most) fractal case $D=0$ (respectively, $D=1$) would correspond to the left-hand side $Re(s)=0$ (respectively, right-hand side $Re(s)=1$) of the critical strip.\,Furthermore, the \emph{mid-fractal case} $D=\frac{1}{2}$ would correspond to the \emph{critical line} $Re(s)=\frac{1}{2}$, along which all of the nontrivial (or critical) zeroes of $\zeta$ are supposed to lie.\\

\hspace*{3mm}The above intuition was both used and justified in the work of \cite{LaMa1, LaMa2}.\,In particular, a key result of \cite{LaMa2} was proved by assuming that $\omega=D+i\tau$ ($\tau>0$) is a zero of $\zeta$ (which implies that $\zeta(\overline{\omega})=0$, where $\overline{\omega}=D-i\tau$), then showing that it follows that $W(x)-N_{\nu}(x)$ is asymptotically proportional to $\mathcal{M}(\mathcal{L})x^{D}$, and finally constructing a fractal string $\mathcal{L}$ of dimension $D$ which is not Minkowski measurable (in light of the above characterization of Minkowski measurability from \cite{LaPo2}).\footnote{Indeed, according to the construction of \cite{LaMa2}, we have $N_{\mathcal{L}}(x)\sim x^{D}+\beta(x^{\omega}+x^{\overline{\omega}})=x^{D}(1+2\beta cos(\tau\log x))$, for some $\beta>0$ small enough.}\,This fractal string provides a counter-example to the inverse spectral problem considered in \S4.1 (recall that we have assumed here that $\zeta(\omega)=0$, with $Re(\omega)=D$), under the above assumption that $\zeta(s)$ has at least one zero along the vertical line $Re(s)=D$.\,In other words, heuristically, the imaginary part $\tau$ of the `complex dimension' $\omega$ gives rise to geometric oscillations (of leading order $D$), thereby showing that $\mathcal{L}$ is not Minkowski measurable (in light of the Minkowski measurability criterion of \cite{LaPo2}), but the spectral oscillations (also of order $D$) that should be associated with $\tau$ are `killed' because $\zeta(\omega)=\zeta(\overline{\omega})=0$.\,We note that in the language of the theory of complex (fractal) dimensions since then developed in \cite{La-vF1, La-vF2, La-vF3}, the fractal string $\mathcal{L}$ constructed in \cite{LaMa2} has precisely for set of complex dimensions
\begin{equation}
\mathcal{D}_{\mathcal{L}}=\{D, \omega, \overline{\omega}\},\mbox{\quad where $\omega=D+i\tau$}.
\end{equation}
Moreover, the explicit formulas from \cite{La-vF2, La-vF3} (see Theorem \ref{Thm:2} and especially, its consequence at the spectral level, Equation (\ref{Eq:19})) can be used in order to obtain a streamlined proof of the fact that the spectral oscillations of $\mathcal{L}$ disappear in this case, because $\zeta(\omega)=\zeta(\overline{\omega})=0$ and $\omega, \overline{\omega}$ are simple poles of $\zeta_{\mathcal{L}}$; see [\textbf{La-vF3}, Ch.\,9].

\begin{remark}
$($\emph{The higher-dimensional case}$.)$\,The Weyl--Berry conjecture \\\cite{Berr1, Berr2} for the vibrations of fractal drums was partially resolved in \cite{La1} in the case of drums with fractal boundaries $($in any dimension $N\geq1$$)$.\,See also the important earlier work of J.\,Brossard and R.\,Carmona \emph{[\textbf{BroCa}]} where was provided a counter-example to the original conjecture $($expressed in terms of the Hausdorff instead of the Minkowski dimension of the boundary$)$ and a corresponding, but weaker, error estimate was obtained for the asymptotics of the trace of the heat semigroup $($or `partition function'$)$, in the special case of the Dirichlet Laplacian.\\
\hspace*{3mm}Accordingly, it was shown by the second author in \cite{La1} that if $\Omega\subset \mathbb{R}^{N}$ $($$N\geq1$$)$ is an arbitrary bounded open set with $($inner$)$ Minkowski dimension $D\in (N-1,N)$ and of finite upper Minkowski content $($i.e., $\mathcal{M}^{*}=\mathcal{M}^{*}_{D}(\partial \Omega)<\infty$$)$,\footnote{We always have that $D\in[N-1,N]$, a statement which reduces to the familiar condition $D\in [0,1]$ in the case of an ordinary fractal string $($i.e., $N=1$$)$; see \cite{La1}, where the cases $D=N-1, N$ and $N-\frac{1}{2}$ are respectively referred to as the least, most and mid-fractal cases.\,Furthermore, we note that in $\mathbb{R}^{N}$, the Minkowski dimension and (upper, lower) Minkowski content are defined exactly as in Definition \ref{Def:Minmeas} above, except with $1-D$ and $1-\alpha$ replaced by $N-D$ and $N-\alpha$, in Equation 
(\ref{Eq:Mm}) and Equation (\ref{Eq:upMink}), respectively.} we have the following remainder estimate for the Dirichlet Laplacian on $\Omega$ $($interpreted either variationally or distributionally$)$$:$
\begin{equation}
N_{\nu}(x)=W(x)+O(x^{D})\mbox{\quad \emph{as} \,$x\to +\infty$,}\label{Eq:Weylerr}
\end{equation}
where $W(x):=(2\pi)^{-N}\mathcal{B}_{N}vol_{N}(\Omega)x^{N}$ is the Weyl $($or leading$)$ term\footnote{When $N=1$, it reduces to the Weyl term given in Equation (\ref{Eq:Weylter}).} with $vol_{N}(\Omega)$ and $\mathcal{B}_{N}$ respectively denoting the $N$-dimensional volume of $\Omega$ and the closed unit ball of $\mathbb{R}^{N}$.\footnote{In the least fractal case where $D=N-1$, the error term on the right-hand side of Equation (\ref{Eq:Weylerr}) should be replaced with $O(x^{N-1}\log x)$.}\,Furthermore, in \cite{La1}, the error estimate in Equation $($\ref{Eq:Weylerr}$)$ is shown to be sharp in every possible dimension $D\in (N-1, N)$.\,Moreover, analogous results are obtained in \cite{La1} for the Neumann Laplacian $($under suitable assumptions on $\partial\Omega$$)$ as well as for positive elliptic operators of order $2m$ $($$m\in \mathbb{N}$, $m\geq1$$)$ and with possibly variable coefficients.\\
\hspace*{3mm}For further discussion of the Weyl--Berry conjecture $($and its later modifications, beginning with \cite{La1} and \cite{LaPo3}$)$ or its physical motivations, we refer, for example, to \emph{[\textbf{Berr1, Berr2, BroCa, La1, La2, La3, La4, LaPo2, LaPo3, FlVa, Ger, GerScm1, GerScm2, HeLa, MolVai, vB-Gi, HamLa}]}, along with \emph{[\textbf{La-vF2}, \S12.5]} and the relevant references therein.\,$($See also, for instance, \emph{[\textbf{Berr1, Berr2, FukSh, La3, KiLa1, Ham1, Ham2, KiLa2, Ki, Sab1, Sab2, Sab3, Str, Tep1, Tep2}]} and the relevant references therein for the case of a drum with a fractal membrane rather than with a fractal boundary.$)$
\end{remark}

\section{Quasi-Invertibility and Almost Invertibility of the Spectral Operator} 
\hspace*{3mm}In order to study the invertibility of the spectral operator $\mathfrak{a}_{c}$, we first introduce two new families of truncated operators: the truncated infinitesimal shifts $\partial_{c}^{(T)}$ and the truncated spectral operators $\mathfrak{a}_{c}^{(T)}$.\, These are the key mathematical objects behind the existence of two notions of invertibility of the spectral operator $\mathfrak{a}_{c}$ which were introduced and studied in \cite{HerLa1}, namely, \emph{quasi-invertibility} and \emph{almost invertibility} (see Definitions \ref{Def:qi} and \ref{Def:ai} below).\,We show in \cite{HerLa1} that these two notions of invertibility play a key role in unraveling the precise relation between the existence of a suitable `inverse' for the spectral operator and the inverse spectral problem for fractal strings studied in \cite{LaMa1, LaMa2} (as well as later on, in \cite{La-vF1, La-vF2, La-vF3}, via the explicit formulas) and discussed in \S4.1.

\subsection{The truncated operators $\partial^{(T)}_{c}$ and $\mathfrak{a}^{(T)}_{c}$.}
\hspace*{3mm}In order to define the notion of quasi-invertibility, we first introduce the \emph{truncated infinitesimal shift} $\partial^{(T)}$ and the \emph{truncated spectral operator} $\mathfrak{a}^{(T)}$.\,As is stated in Corollary \ref{CorTh} and Theorem \ref{Thm:SpVc} (see Appendix B) or follows equivalently from Theorems \ref{Thm:part} and \ref{Thm:spp}, the infinitesimal shift $\partial=\partial_{c}$ is given by
\begin{equation}
\partial_{c}=c+iV,
\end{equation}
\text
where $V=V_{c}$ is an unbounded self-adjoint operator on $\mathbb{H}_{c}$ with spectrum $\sigma(V)=\mathbb{R}$.\,Thus, given $T\geq 0$, we define the \emph{truncated infinitesimal shift} as follows:
\begin{equation}
A^{(T)}=\partial^{(T)}:=c+iV^{(T)},
\end{equation}
where
\begin{equation}
V^{(T)}:=\phi^{(T)}(V)\notag\\
\end{equation}
\text
and $\phi^{(T)}$ is a suitable (i.e., $T$-admissible) cut-off function (so that we have, in particular, $\sigma(A^{(T)})=c+i[-T,T]$).
\begin{remark}
More precisely, $\phi^{(T)}$ is any $T$-\emph{admissible cut-off function}, defined as follows$:$ when $c\ne 1$, $\phi^{(T)}$ is a continuous function defined on $\mathbb{R}$ and the closure of its range is equal to $[-T,T]$.\,Furthermore, when $c=1$, $\phi^{(T)}$ is meromorphic in an open neighborhood of $\mathbb{R}$ in $\mathbb{C}$ and the closure of the range of its restriction to $\mathbb{R}$ is equal to $[-T,T]$; in this case, one views $\phi^{(T)}$ as a continuous function with values in the Riemann sphere $\widetilde{\mathbb{C}}:=\mathbb{C}\cup\{\infty\}$.\,$($For example, we may take
$\phi^{(T)}(s)=\frac{T}{\pi}\tan^{-1}(s)$, initially defined for $s\in \mathbb{R}$.$)$\,One then uses the measurable functional calculus for unbounded normal operators, along with the corresponding continuous $($$c\ne 1$$)$ or meromorphic $($$c=1$$)$ version of the spectral mapping theorem $($see the relevant appendix in \cite{HerLa1} and Remark \ref{Rk:spt} above$)$ in order to define both $\partial^{(T)}$ and $\mathfrak{a}^{(T)}$ and calculate their spectra.
\end{remark}

\hspace*{3mm}Similarly, in light of the definition of the (standard) spectral operator $\mathfrak{a}=\mathfrak{a}_{c}$ given in Equation (\ref{Def:Spop}), the \emph{truncated spectral operator} $\mathfrak{a}^{(T)}=\mathfrak{a}^{(T)}_{c}$ is defined by
\begin{equation}
\mathfrak{a}_{c}^{(T)}:=\zeta\left(\partial^{(T)}\right).
\end{equation}

\hspace*{3mm}Note that the above construction can be generalized as follows: \\

Given $0\leq T_{0}\leq T$, one can define a $(T_{0},T)$-\emph{admissible cut-off function} $\phi^{(T_{0},T)}$ exactly as above, except with $[-T,T]$ replaced with $\{\tau\in\mathbb{R}:\,T_{0}\leq |\tau| \leq T\}$.\\

Correspondingly, one can define $V^{(T_{0},T)}=\phi^{(T_{0},T)}(V)$,\\
\begin{equation}
A^{(T_{0},T)}=\partial^{(T_{0},T)}:=c+iV^{(T_{0},T)}\label{Trunc1}
\end{equation}
and
\begin{equation}
\mathfrak{a}_{c}^{(T_{0},T)}=\zeta(\partial^{(T_{0},T)}),\label{Trunc2}
\end{equation}
where $\partial^{(T_{0},T)}$ is the $(T_{0},T)$-\emph{infinitesimal shift} and $\mathfrak{a}^{(T_{0},T)}$ is the $(T_{0},T)$-\emph{truncated spectral operator}.

\begin{remark}
Note that when we let $T_{0}=0$ in Equations \emph{(\ref{Trunc1})} and \emph{(\ref{Trunc2})}, we recover $A^{(T)}$ and $\mathfrak{a}_{c}^{(T)}$; i.e., $A^{(T)}=A^{(0,T)}$ and $\mathfrak{a}_{c}=\mathfrak{a}_{c}^{(0,T)}$.
\end{remark}

\hspace*{3mm}Finally, we introduce the notions of \emph{quasi-invertibility}  and \emph{almost invertibility} of $\mathfrak{a}=\mathfrak{a}_{c}$ as follows (the standard notion of invertibility of an operator is recalled in Remark \ref{Rk::Invunb} below):

\begin{definition}\label{Def:qi}
The spectral operator $\mathfrak{a}$ is \emph{quasi-invertible} if its truncation $\mathfrak{a}^{(T)}$ is invertible for all $T>0$.
\end{definition}

\begin{definition}\label{Def:ai}
Similarly, $\mathfrak{a}$ is \emph{almost invertible} if for some $T_{0}\geq 0$, its truncation $\mathfrak{a}^{(T_{0},T)}$ is invertible for all $T>T_{0}$.
\end{definition}

Note that in the definition of \textquotedblleft almost invertibility\textquotedblright, $T_{0}$ is allowed to depend on the parameter $c$.\,Furthermore, observe that quasi-invertiblity implies almost invertibility.\\
\begin{remark}\label{Rk::Invunb}
Recall that a $($possibly unbounded$)$ densely defined linear operator $A:\,D(A)\subset \mathbb{H}\to \mathbb{H}$ on a Hilbert space $\mathbb{H}$, where $D(A)$ is the domain of $A$, is said to be \emph{invertible} if it is invertible in the set theoretic sense and if its inverse is bounded.\footnote{If, in addition, $A$ is closed (which will be the case of all of the operators considered here, including $\partial_{c}$, $\mathfrak{a}_{c}$ and its truncations $\mathfrak{a}^{(T)}_{c}$ and $\mathfrak{a}^{(T, T_{0})}_{c}$), the inverse operator is automatically bounded (by the closed graph theorem); see, e.g., \cite{DunSch, Kat, Ru}.}In other words, there exists a bounded linear operator $B:\,\mathbb{H}\to\mathbb{H}$ with range $D(A)$ and defined on all of \,$\mathbb{H}$, such that $ABu=u$ for all $u\in\mathbb{H}$ and $BAv=v$, for all $v\in D(A)$.\,Furthermore, note that according to the definition of the spectrum $\sigma(A)$ of $A$, the linear operator $A$ is invertible if and only if $0\notin \sigma(A)$ $($See, e.g., \cite{DunSch, Kat, ReSi, Ru, Sc}.$)$ 
\end{remark}

\subsection{The spectra of $\partial_{c}^{(T)}$ and $\mathfrak{a}^{(T)}_{c}$.}

\hspace*{3mm}The spectra of $A^{(T)}$ and $\mathfrak{a}^{T}$ are now respectively determined as follows:

\begin{theorem}\emph{\cite{HerLa1}}\,For all $T>0$, $A^{(T)}$ is a bounded normal linear operator whose spectrum is given by \label{Thm:spcqt}
\begin{equation}
\sigma(A^{(T)})=\{c+i\tau:\,|\tau|\leq T\}.\label{Eq:trpr}
\end{equation}
\end{theorem}

\begin{theorem}\label{Thm:spbounnboun}\emph{\cite{HerLa1}}\,For all $T>0$, $\mathfrak{a}_{c}^{(T)}$ is a bounded normal linear operator \footnote{More precisely, only when $c=1$, which corresponds to the pole of $\zeta(s)$ at $s=1$, $\mathfrak{a}^{(T)}$ is not bounded (since $\zeta(1)=\infty\in \widetilde{\mathbb{C}}$) and hence, Equation (\ref{Eq:trsp}) must then be interpreted as an equality in $\widetilde{\mathbb{C}}$, with $\zeta$ viewed as a $\widetilde{\mathbb{C}}$-valued (continuous) function.\,(See Remark \ref{Rk:spt}.)} whose spectrum is given by
\begin{equation}
\sigma(\mathfrak{a}_{c}^{(T)})=\{\zeta(c+i\tau):\,|\tau|\leq T \}.\label{Eq:trsp}
\end{equation}
\end{theorem}

\hspace*{3mm}More generally, given $0\leq T_{0}\leq T$, the exact counterpart of Theorem \ref{Thm:spcqt} and Theorem \ref{Thm:spbounnboun} holds for $A^{(T_{0},T)}$ and $\mathfrak{a}^{(T_{0},T)}=\mathfrak{a}_{c}^{(T_{0},T)}$, except with $|\tau|\leq T$ replaced with $T_{0}\leq |\tau|\leq T$.\\

\hspace*{3mm}Our next result provides a necessary and sufficient condition for the invertibility of the truncated spectral operator:\footnote{Recall from the end of Remark \ref{Rk::Invunb} that by definition of the spectrum, $\mathfrak{a}_{c}^{(T)}$ is invertible if and only if $0\notin \sigma(\mathfrak{a}_{c}^{(T)})$.}\\

\begin{corollary}\emph{\cite{HerLa1}}\,Assume that $c\geq 0$.\,Then, the truncated spectral operator $\mathfrak{a}^{(T)}$ is invertible if and only if $\zeta$ does not have any zeroes on the vertical line segment $\{s\in \mathbb{C}:\,Re(s)=c,\,|Im(s)|\leq T\}$.\label{Cor:trun}
\end{corollary}

Naturally, given $0\leq T_{0}\leq T$, the same result as in Corollary \ref{Cor:trun} is true for $\mathfrak{a}^{(T_{0},T)}$ provided $|Im(s)|\leq T$ is replaced with $T_{0}\leq |Im(s)|\leq T$.\\

\subsection{Quasi-invertibility of $\mathfrak{a}_{c}$, almost invertibility and Riemann zeroes.}

\hspace*{3mm}Next, we deduce from the above results necessary and sufficient conditions ensuring the quasi-invertibility or the almost invertibility of $\mathfrak{a}=\mathfrak{a}_{c}$.\,Such conditions turn out to be directly related to the location of the critical zeroes of the Riemann zeta function.

\begin{theorem}\emph{\cite{HerLa1}}\,Assume that $c\geq 0$.\,Then, the spectral operator $\mathfrak{a_{c}}=\zeta(\partial_{c})$ is quasi-invertible if and only if the Riemann zeta function does not vanish on the vertical line $\{s\in \mathbb{C}:\,Re(s)=c\}$.\label{Thm:qinv}
\end{theorem} 
 \hspace*{3mm}We now state the exact counterpart of Theorem \ref{Thm:qinv} for the almost invertibility (rather than the quasi-invertibility) of $\mathfrak{a}=\mathfrak{a}_{c}$.
 
\begin{theorem}\emph{[\textbf{HerLa1}]}\,Assume that $c\geq0$.\,Then, $\mathfrak{a}_{c}$ is almost invertible if and only if all but $($at most$)$ finitely many zeroes of $\zeta$ are off the vertival line $Re(s)=c$.\label{Thm:ai}
\end{theorem} 

\begin{remark}
In light of Definition \ref{Def:qi}, Theorem \ref{Thm:qinv} follows from Corollary \ref{Cor:trun}.\,Similarly, in light of Definition \ref{Def:ai}, Theorem \ref{Thm:ai} follows from the counterpart $($or really, the extension$)$ of Corollary \ref{Cor:trun} mentioned in the comment following that corollary.\,Moreover, it is worth pointing out that the definitions of the Hilbert space $\mathbb{H}_{c}$ as well as of the spectral $\mathfrak{a}=\mathfrak{a}_{c}$ and of its truncations $\mathfrak{a}^{(T)}_{c}$ and $\mathfrak{a}_{c}^{(T, T_{0})}$ given in \cite{HerLa1} make sense for any $c\in\mathbb{R}$.\,Accordingly, all of the results stated in \S5 have an appropriate counterpart for any $c\in\mathbb{R}$ provided we take into account the trivial zeroes of $\zeta(s)$, located at $s=-2n$, for $n=1, 2,...$\,For the simplicity of exposition, we will not further discuss this issue here.\,$($See Appendix B, however.$)$
\end{remark}

\section{Spectral Reformulations of the Riemann Hypothesis and of Almost RH}
\subsection{Quasi-invertibility of $\mathfrak{a}_{c}$ and spectral reformulation of RH}

\hspace*{3mm}In this subsection, we first deduce from our earlier results in \S5.3 (specifically, from Theorem \ref{Thm:qinv}) a \emph{spectral reformulation of the Riemann hypothesis} (RH, see Theorem \ref{Thm:qual} below), expressed in terms of the quasi-invertibility of the spectral operator $\mathfrak{a}=\mathfrak{a}_{c}$.\,From a functional analytic and operator theoretic point of view, this reformulation sheds new light on, and further extends, the work of the second author and H.\,Maier \cite{LaMa2} in their study of \emph{the inverse spectral problem for vibrating fractal strings}.\,(See \S4.1 above for a brief description of this inverse problem and of the main results of \cite{LaMa2}.)\,This result also sheds new light on the reinterpretation and further extensions of the work of \cite{LaMa2} obtained in [\textbf{La-vF3}, Ch.\,9] in terms of a rigorously formulated theory of complex dimensions and the corresponding explicit formulas.\,(Recall from [\textbf{La-vF3}, \S6.3.1] as well as from \S2.2 and \S3.1 above that the heuristic spectral operator $\eta\mapsto\nu \notag$ can be understood in terms of the explicit formulas of \cite{La-vF2, La-vF3} expressed in terms of the geometric and spectral complex dimensions of generalized fractal strings; see Equation (\ref{Spoplev1}), along with Equations (\ref{Eq:18})and (\ref{Eq:19}).)\,In particular, Theorem \ref{Thm:qual} below enables us to give a precise mathematical meaning in this context to the notion of invertibility of the spectral operator, as discussed semi-heuristically in [\textbf{La-vF3}, Cor.\,9.6].\,Indeed, here, the proper notion of invertibility of $\mathfrak{a}$ is that of quasi-invertibility.

\begin{theorem}\emph{\cite{HerLa1}}\label{Thm:qual}
The spectral operator $\mathfrak{a}=\mathfrak{a}_{c}$ is quasi-invertible for all $c\in(0,1)-\frac{1}{2}$ $($or equivalently, for all $c\in (\frac{1}{2},1)$$)$ if and only if the Riemann hypothesis is true.
\end{theorem}

\begin{remark}\label{Rk:globzet}
The fact that the dimensional parameter $c$ may equivalently be assumed to lie in $(0,\frac{1}{2})$, $(\frac{1}{2},1)$ or all of $(0,1)-\frac{1}{2}$ follows from the functional equation for the Riemann zeta function, which connects $\zeta(s)$ and $\zeta(1-s)$; see Equations \emph{(\ref{Eq:fE})} and \emph{(\ref{Eq:CZ})} in Appendix A.\,$($An entirely analogous comment can be made about Theorem \ref{Thm:ainv} below.$)$ 
\end{remark}

\subsection{Almost invertibility of $\mathfrak{a}_{c}$ and spectral reformulation of \textquotedblleft almost RH\textquotedblright.}
\hspace*{3mm}We next deduce from the results of \S5.3 (specifically, from Theorem \ref{Thm:ai}) a new statement concerning $\zeta$, to which we refer to as a spectral reformulation of the \emph{almost Riemann hypothesis} (almost RH, in short).
\begin{theorem}\emph{\cite{HerLa1}}\label{Thm:ainv}
The spectral operator $\mathfrak{a}=\mathfrak{a}_{c}$ is almost invertible for all $c\in (\frac{1}{2},1)$ if and only if the Riemann hypothesis $($RH$)$ is \textquotedblleft almost true\textquotedblright  $($i.e., on every vertical line $Re(s)=c$, with $c>\frac{1}{2}$, there are at most finitely many exceptions to RH$)$.
\end{theorem}

\begin{remark}
Theorem \ref{Thm:ainv} $($as well as Theorem \ref{Thm:ai}, of which it is a corollary$)$ does not have any counterpart in the results of \cite{LaMa1, LaMa2} and \cite{La-vF2, La-vF3} or, to our knowledge, in the existing reformulations of the Riemann hypothesis and of its many variants.\,Furthermore, recall that $\zeta(s)\ne0$ for $Re(s)\geq1$ $($for $Re(s)=1$, this is Hadamard's theorem$)$; see Appendix A.\,This fact explains why we wrote $c>\frac{1}{2}$ instead of $c\in(\frac{1}{2},1)$ in the latter part of Theorem \ref{Thm:ainv}.
\end{remark}

\begin{remark}\label{Rk:arith}
\hspace*{3mm}Theorem \ref{Thm:qual} and Theorem \ref{Thm:ainv} have natural counterparts for a very large class of arithmetic zeta functions $($or $L$-functions$)$, thereby yielding a new operator theoretic and spectral reformulation of the generalized Riemann hypothesis $($GRH$)$ and of the \textquotedblleft almost GRH\textquotedblright, respectively.\,Naturally, the corresponding generalized spectral operator $\mathfrak{a}_{L,\,c}$ would then be defined by $\mathfrak{a}_{L,\,c}=L(\partial_{c})$, where $L=L(s)$ is the $L$-function under investigation; see \S7.1 and \S7.3 below.
\end{remark}

\begin{remark}\label{RK:qiai}
Note that according to our previous results and definitions, the invertibility of the spectral operator $\mathfrak{a}$ implies its quasi-invertibility, which in turn implies its almost invertibility.
\end{remark}

\hspace*{3mm}In light of Remark \ref{RK:qiai}, we deduce the following corollary from Theorem \ref{Thm:ainv} and Hardy's theorem according to which $\zeta$ has infinitely many zeroes on the critical line $Re(s)=\frac{1}{2}$ (see, e.g., \cite{Tit}).

\begin{corollary}\emph{\cite{HerLa1}}\,For $c=\frac{1}{2}$, the spectral operator $\mathfrak{a}$ is not almost $($and thus, not quasi-$)$ invertible.
\end{corollary}

\subsection{Invertibility of the spectral operator and phase transitions.}

\hspace*{3mm}We have discussed in \S5.3, \S6.1 and \S6.2 various characterizations of the quasi-invertibility or of the almost invertibility of the spectral operator $\mathfrak{a}=\mathfrak{a}_{c}$, either for a given $c\geq0$ in \S5.3, or else for all $c\in(\frac{1}{2},1)$ (or equivalently, for all $c\in(0,1)-\frac{1}{2}$), in \S6.1 and \S6.2, respectively.\,In the present subsection, however, we very briefly discuss the invertibility of $\mathfrak{a}$, in the standard sense of closed (possibly unbounded) operators recalled in Remark \ref{Rk::Invunb}.\,As it turns out, one has to distinguish three main cases: $c>1$, $\frac{1}{2}<c<1$ and $0<c<\frac{1}{2}$.\,Although, there are several very interesting new features that are quite different from those encountered in \S5.3, \S6.1 and \S6.2, our discussion of the (standard) invertibility of $\mathfrak{a}$ will be rather succinct  because it is not the main object of the present paper.\,A detailed discussion can be found in \cite{HerLa1} and a survey of this topic is provided in \cite{HerLa3}.\\

\hspace*{3mm}Recall from the end of Remark \ref{Rk::Invunb} that, by definition of the spectrum $\sigma(\mathfrak{a})$ of $\mathfrak{a}$, the operator $\mathfrak{a}$ is invertible if and only if $0\notin \sigma(\mathfrak{a})$.\,We therefore deduce from Theorem \ref{Thm:ssop} (the characterization of the spectrum of $\mathfrak{a}$) the following \emph{invertibility criterion} for $\mathfrak{a}$.

\begin{theorem}\emph{\cite{HerLa1}}\,Assume that $c\geq 0$.\,Then, the spectral operator $\mathfrak{a}$ is invertible if and only if $0\notin cl(\{\zeta(s):\,Re(s)=c\})$.\footnote{That is, if and only if $\zeta$ does not have any zeroes on the vertical line $L_{c}:=\{Re(s)=c\}$ and there is no infinite sequence $\{s_{n}\}_{n=1}^{\infty}$ of distinct points of $L_{c}$ such that $\zeta(s_{n})\to 0$ as $n\to \infty$.}\label{Thm:inv}
\end{theorem}

\hspace*{3mm}Next, we will explore some of the consequences of Theorem \ref{Thm:quasinv} in light of the universality of $\zeta(s)$ in the right critical strip $\{\frac{1}{2}<Re(s)<1\}$ (see part (2) of Theorem \ref{Thm:quasinv}) and conditionally (i.e., under RH) of the non-universality of $\zeta(s)$ in the left critical strip $\{0<Re(s)<\frac{1}{2}\}$\,(see part $(3)$ of Theorem \ref{Thm:quasinv}, which makes use of the work of R. Garunk\v{s}tis and J. Steuding in \cite{GarSt}).\\

\begin{theorem}\emph{\cite{HerLa1, HerLa3}}\,Assume that $c>0$.\,Then$:$ \label{Thm:quasinv}
\begin{enumerate}
\item For $c>1$, $\mathfrak{a}$ is invertible\footnote{In light of Remark \ref{RK:qiai}, it follows that $\mathfrak{a}$ is also quasi-\,(and hence, almost) invertible.} and bounded; its spectrum is a compact subset of $\mathbb{C}$ avoiding the origin.
\item For $c\in(\frac{1}{2},1)$, $\mathfrak{a}$ is not invertible and in fact, $\sigma(\mathfrak{a})=\mathbb{C}$.\,In particular, $\mathfrak{a}$ is unbounded.
\item For $c\in(0,\frac{1}{2})$, $\mathfrak{a}$ is also unbounded $($i.e., $\sigma(\mathfrak{a})$ is unbounded$)$ and, assuming the Riemann hypothesis $($i.e., conditionally$)$, $\mathfrak{a}$ is not invertible $($i.e., $0\notin \sigma(\mathfrak{a})$$)$.\footnote{It is not known whether the conclusion of part $(3)$ is true unconditionally or can be drawn under a weaker hypothesis than RH; see \cite{HerLa1, HerLa3}.\,(See also \cite{GarSt}.)}
\end{enumerate}
\end{theorem}

\hspace*{3mm}As was alluded to above, Theorem \ref{Thm:quasinv} exhibits two different types of (mathematical) phase transitions, one occurring at $c=1$, and conditionally, another one occurring at $c=\frac{1}{2}$.\,These `phase transitions' correspond to both the nature (or the shape) of the spectrum, the boundedness of $\mathfrak{a}$,\footnote{In fact, it follows from Theorems \ref{Thm:inv} and \ref{Thm:quasinv} that $\mathfrak{a}$ is bounded for $c>1$ and unbounded for $0<c\leq 1$.} and the invertibility of $\mathfrak{a}$.\,The possible physical origins and interpretations of these phase transitions are discussed in \cite{HerLa1} and \cite{HerLa3}.\\

\hspace*{3mm}We note that the spectral reformulation of the Riemann hypothesis provided in \S6.1 (and that of \textquotedblleft almost RH\textquotedblright provided in \S6.2) is associated with yet another (mathematical) phase transition, occurring this time only at $c=\frac{1}{2}$ (which corresponds, of course, to both the \emph{mid-fractal dimension} $D=\frac{1}{2}$ and the \emph{critical line} $Re(s)=\frac{1}{2}$).\,The same comment can be made about the earlier reformulations of RH obtained in \cite{LaMa1, LaMa2} and later on, in \cite{La-vF2, La-vF3}; see \emph{loc.\,cit}. and \cite{La3}.\\

\begin{remark}
We have seen above that the issue of universality of the Riemann zeta function $($and, more generally, of other $L$-functions, see \cite{St}$)$ plays an important role in aspects of the present theory.\footnote{We refer the interested reader to the books \cite{KarVo}, \cite{Lau} and \cite{St} for an exposition of the theory of universality, originating $($in the 1920s and in the 1970s, respectively$)$ with the beautiful Bohr--Landau and Voronin theorems.\,We simply mention here that roughly speaking, \textquotedblleft universality\textquotedblright $($in this context$)$ means that any non-vanishing holomorphic function (in a suitable compact subset of $\{\frac{1}{2}<Re(s)<1\}$) can be uniformly approximated by vertical translates of $\zeta$ (or of the given $L$-function under study).}\,This topic is explored in \cite{HerLa1} and in \cite{HerLa4} where the spectral operator $\mathfrak{a}=\zeta(\partial)$, viewed as a suitable \emph{quantization} of the Riemann zeta function, is shown to be \textquotedblleft universal\textquotedblright $($in an appropriate sense$)$ among all non-vanishing holomorphic functions\footnote{restricted to a suitable compact subset of the right critical strip $\{\frac{1}{2}<Re(s)<1\}$.} of the infinitesimal shift $\partial=\partial_{c}$ $($which now plays the role of the complex variable $s$ in the classic theory of universality$)$.\footnote{The actual formulation of the universality of $\mathfrak{a}$ is a little more complicated and involves the family of truncated spectral operators $\mathfrak{a}^{(T)}$.} 
\end{remark}

\begin{center}
\section{Concluding Comments}
\end{center}
\hspace*{3mm}The functional analytic framework which was provided in \cite{HerLa1} was crucial to give a precise mathematical meaning to the heuristic definition of the spectral operator given in [\textbf{La-vF3}, \S6.3].\,Indeed, it enabled us to rigorously define and study the infinitesimal shift $\partial_{c}$, the shift (or translation) semigroup $e^{t\partial_{c}}$, the spectral operator $\mathfrak{a}_{c}=\zeta(\partial_{c})$ and its appropriate truncations $\mathfrak{a}^{(T)}_{c}$ (and $\mathfrak{a}^{(T_{0},\,T)}_{c}$),\footnote{We note that the notion of `truncated spectral operator' does not appear in \cite{La-vF3}.\,In fact, we were led to introducing it in \cite{HerLa1} in order to find the appropriate notion of invertibility (namely, quasi-invertibility) necessary to obtain this reformulation.} determine their spectra and thus, obtain the spectral reformulation of the Riemann hypothesis (RH) in Theorem \ref{Thm:qual} while investigating the invertibility of the spectral operator. As a result, an extension of, and a new operator theoretic perspective on, the earlier work in \cite{LaMa2} were obtained.\\

\subsection{Extension to arithmetic zeta functions.}
\hspace*{3mm}As was alluded to earlier, the criteria provided in Theorems \ref{Thm:qinv}, \ref{Thm:ai}, \ref{Thm:inv} and \ref{Thm:quasinv} clearly extend in a natural manner to the spectral operators associated to a large class of  arithmetic zeta functions $($or L-functions$)$.\,The same can be said of most of the results of \cite{HerLa1, HerLa2, HerLa3, HerLa4} discussed in this survey.\\

\subsection{Operator-valued Euler products.}
\hspace*{3mm}Furthermore, one can show $($see \cite{HerLa2}$)$ that for $c>1$, $\mathfrak{a}_{c}$ belongs to $\mathcal{B}(\mathbb{H}_{c})$ and is given by the following \emph{operator-valued Euler product expansion} for $\mathfrak{a}=\mathfrak{a}_{c}$:

\begin{equation}
\mathfrak{a}_{c}=\zeta(\partial)=\prod_{p\in\mathcal{P}}(1-p^{-\partial})^{-1},\label{Eq:EPrd}
\end{equation}
\text
where $\partial=\partial_{c}$ and the convergence holds in the Banach algebra $\mathcal{B}(\mathbb{H}_{c})$ of bounded linear operators on $\mathbb{H}_{c}$.\,Moreover, still for $c>1$, we have that $||\mathfrak{a}_{c}||\leq \zeta(c)$ and $\mathfrak{a}$ is invertible with $($bounded$)$ inverse given by

\begin{equation}
\mathfrak{a}_{c}^{-1}=\frac{1}{\zeta}(\partial)=\sum_{n=1}^{\infty}\mu(n)n^{-\partial},\label{Eq:3.10}
\end{equation}
\text
where the equality holds in $\mathcal{B}(\mathbb{H}_{c})$ and $\mu(n)$ is the classic M\"obius function defined by $\mu(n)=(-1)^{q}$ if $n\in \mathbb{N}$ is a product of $q$ distinct primes, and $\mu(n)=0$, otherwise.\,$($Compare Equations (\ref{Eq:3.10}) and (\ref{Eq:2.23}).\,Also, recall that for $c>1$, Equation (\ref{Eq:2.23}) was rigorously justified by Theorem \ref{Thm:11}.$)$\,In addition, it was conjectured in [\textbf{La-vF3}, \S6.3.2] that the above Euler product in Equation (\ref{Eq:EPrd}) also converges $($in a suitable sense$)$ inside the critical strip, that is, for $0<c<1$.\,This conjecture is addressed in \cite{HerLa2}.\\

\subsection{Global spectral operator.}
\hspace*{3mm}Another interesting problem consists in considering and studying the\emph{ global spectral operator} $\mathcal{A}_{c}:=\xi(\partial_{c})$, where $\xi$ is the global (or completed) Riemann zeta function given in Equation (\ref{Eq:CZ}) of Appendix A and satisfies the functional equation (\ref{Eq:fE}): $\xi(s)=\xi(1-s)$.\,Due to the perfect symmetry of the functional equation, this operator has some appealing properties, particularly for $c=\frac{1}{2}$.\,In particular, an operator-valued functional equation connecting $\mathcal{A}_{c}$ and $\mathcal{A}_{1-c}$ can be obtained (see \cite{HerLa1}).\,Naturally, in the spirit of \S7.1, an analogous problem can be investigated for \emph{generalized global spectral operators} associated with global (or completed) $L$-functions.\\

\subsection{Towards a quantization of number theory.}
\hspace*{3mm}In closing, we note that our study of the spectral operator provides a `natural quantization' of several identities in analytic number theory.\,For instance, as was briefly discussed at the end of \S6.3, we show in \cite{HerLa1} (see also \cite{HerLa4}) that one can obtain a `quantization' of Voronin's theorem about the universality of the Riemann zeta function which states that any non-vanishing holomorphic function in a compact subset of $\{\frac{1}{2}<Re(s)<1\}$ can be uniformly approximated by imaginary translates of $\zeta=\zeta(s)$.\,In our context, and as a consequence, the `universality of the spectral operator $\mathfrak{a}=\zeta(\partial)$' will imply that the spectral operator can emulate any type of complex behavior.\,As a result, it is \emph{chaotic and fractal} \cite{HerLa1, HerLa4}.\,(Possible connections with various aspects of the research program developed in the book \cite{La5} still need to be explored in this context; see also the work in preparation \cite{La6}.) 

\section{Appendix A:\,Riemann's Explicit Formula}

\hspace*{3mm}In this appendix, we first recall for the non-expert some basic properties of the Riemann zeta function $\zeta$.\,We then briefly discuss Riemann's explicit formula and explain the underlying `duality' between the prime powers and the zeroes of $\zeta$.\,Finally, we point out the analogy between Riemann's explicit formula and the (generalized) explicit distributional formulas of \cite{La-vF2, La-vF3} recalled in Theorem \ref{Thm:2}.\,Indeed, in the latter formulas, the underlying `duality' is now between a generalized fractal string $\eta$ \footnote{Here, we point out, in particular, the special case for which a generalized fractal string $\eta=\sum_{j=1}^{\infty}w_{l_{j}}\delta_{l_{j}^{-1}}$ is viewed (in the distributional sense) as an object encoding the geometry of a standard fractal string $\mathcal{L}=\{l_{j}\}_{j=1}^{\infty}$ with lengths (or scales) $l_{j}$ and multiplicities $w_{l_{j}}$.} and its associated complex dimensions.\\

\hspace*{3mm}We recall that Riemann showed in his celebrated 1858 paper \cite{Rie} on the distribution of prime numbers that 
\begin{equation}\label{Eq:Rmprd}
\zeta(s)=\sum_{n=1}^{\infty}n^{-s}=\prod_{p=1}^{\infty}\frac{1}{1-p^{-s}},\mbox{\quad for $Re(s)>1$}
\end{equation}
and that $\zeta$ has a meromorphic continuation to all of $\mathbb{C}$ with a single (and simple) pole at $s=1$, which satisfies the \emph{functional equation}
\begin{equation}
\xi(s)=\xi(1-s),\mbox{\,}\,s\in \mathbb{C},\label{Eq:fE}
\end{equation}
where
\begin{equation}
\xi(s):=\pi^{-\frac{s}{2}}\Gamma(\frac{s}{2})\zeta(s)\label{Eq:CZ}
\end{equation}
is the \emph{completed} $($or \emph{global}$)$ Riemann zeta function (Here, $\Gamma$ denotes the classic gamma function.)\,Note that the trivial zeros of $\zeta(s)$ at $s=-2n$ for $n=1, 2, 3, ...,$ correspond to the poles of the gamma function $\Gamma(\frac{s}{2})$.\,Riemann also conjectured that the nontrivial $($or \emph{critical}$)$ zeros of $\zeta(s)$ $($i.e., the zeros of $\zeta(s)$ which are located in the critical strip $0<Re(s)<1$$)$ all lie on the \emph{critical line} $Re(s)=\frac{1}{2}$.\,This famous conjecture is known as the \emph{Riemann hypothesis}.\\

\hspace*{3mm}It is well known that the Euler product in Equation (\ref{Eq:Rmprd}) converges absolutely to $\zeta(s)$ for $Re(s)>1$ and also uniformly on any compact subset of the half-plane $Re(s)>1$.\,As a result, $\zeta(s)$ does not have any zeroes for $Re(s)>1$.\,Now, using the `symmetry' expressed by the functional equation (\ref{Eq:fE}), one deduces at once that the Riemann zeta function does not have any other zeroes in the region $Re(s)<0$, except for the `non-critical' (or trivial) zeroes corresponding to the poles of the gamma function $\Gamma(\frac{s}{2})$.\\

\hspace*{3mm} In 1892, Hadamard showed that $\zeta(s)$ does not have any zeroes on the vertical line $Re(s)=1$.\,A few years later, in 1896, his result turned out to be a key step in the proof of the Prime Number Theorem.\,Hence, and again using the symmetry of the functional equation (\ref{Eq:fE}), one can conclude that $\zeta(s)$ does not have any zeroes on the vertical line $Re(s)=0$.\,It follows that the critical strip (i.e., the subset $0<Re(s)<1$) is the region of the complex plane in which the \emph{nontrivial} zeroes of $\zeta(s)$ are located.\,Moreover, we point out the fact that in light of Equation (\ref{Eq:CZ}) and the properties of the meromorphic continuation of $\zeta(s)$, the zeroes of $\xi(s)$ coincide with the critical zeroes of $\zeta(s)$ and these zeroes come in complex conjugate pairs (really, in 4-tuples, due to (\ref{Eq:fE}) and provided they do not lie on the critical line), in the critical strip.\footnote{Naturally, $\xi$ is meromorphic in all of $\mathbb{C}$, with two (simple) poles located at $s=0$ and $s=1$.}\\

\hspace*{3mm}In his same 1858 paper, Riemann obtained an \emph{explicit formula} connecting an expression involving the prime numbers $($for example, the prime number counting function$)$, on the one hand, and the $($trivial and critical$)$ zeroes of the Riemann zeta function $\zeta(s)$, on the other hand.\footnote{We refer, for example, the interested reader to \cite{Edw, Ing, Ivi, Pat, Tit, La5, La-vF2, La-vF3} for more detailed information about Riemann's original explicit formula and its various number theoretic generalizations.}\\

\hspace*{3mm}Consider the counting function $f(x):=\sum_{p^{n}\leq x} \frac{1}{n}$, defined for $x>0$.\,In other words, $f(x)$ is the number of prime powers  $p^{n}$ ($n\in\mathbb{N}$, $n\geq1$) not exceeding $x$, each counted with a weight $\frac{1}{n}$.\,Then, a modern version of Riemann's explicit formula can be stated as follows:
\begin{equation}\label{Eq:Rex}
f(x)=\sum_{p^{n}\leq x}\frac{1}{n}=Li(x)-\sum_{\rho}Li(x^{p})-\int_{x}^{+\infty}\frac{1}{t^{2}-1}\frac{dt}{t\log t}-\log2,
\end{equation}
where $n$ runs through all positive integers, $Li(x):=\int_{0}^{x}\frac{dt}{\log t}$ is the logarithmic integral and $\rho$ runs through all the zeroes of the Riemann zeta function, taken in order of increasing absolute values (and for the critical zeroes, in complex conjugate pairs).\,Note that Equation (\ref{Eq:Rex}) provides a `\emph{duality}' between the integral powers of the primes and the zeroes of zeta.\footnote{Actually, for pedagogical reasons, we do not give here the more complicated Riemann explicit formula in its original form, which was expressed in terms of the standard prime number counting function.}\\

\hspace*{3mm}This duality between the primes $p$ $($or additively, their logarithms $\log p$$)$ and the zeroes $($and the pole$)$ of $\zeta(s)$ has been key to most approaches to the Riemann hypothesis.\,In a similar spirit, the generalization of Riemann's explicit formula obtained in \cite{La-vF2,La-vF3} and discussed earlier in Theorem \ref{Thm:2} connects certain expressions involving a generalized fractal string $\eta$ $($for example, the geometric or the spectral counting function of $\eta$$)$, on the one hand, and the geometric or spectral complex dimensions of $\eta$, on the other hand; that is, the poles of the geometric or spectral zeta function of $\eta$.\footnote{Note that the zeroes and the pole of $\zeta$ (along with their multiplicities$)$ can be recovered from the poles $($and the sign of the residues$)$ of the logarithm derivative $-\frac{\zeta'(s)}{\zeta(s)}$.}

\section{Appendix B:\,The Momentum Operator and Normality of $\partial_{c}$}

\hspace*{3mm}The goal of this appendix is to provide the main steps of a proof of Theorem \ref{Thm:part} and then to explain how to deduce from Theorem \ref{Thm:part} (to be reformulated in Corollary \ref{CorTh} below) the characterization of the spectrum of $\partial_{c}$ obtained in Theorem \ref{Thm:spp} (and to be reformulated in Theorem \ref{Thm:SpVc} below).\,The aforementioned restatements of Theorems \ref{Thm:part} and \ref{Thm:spp} will be expressed in terms of the $c$-\emph{momentum operator} $V_{c}$, which we define next.\\

\hspace*{3mm}We recall that the infinitesimal shift $\partial_{c}$ was studied in detail in \cite{HerLa1} and that some of its fundamental properties were presented in \S3.3 above.\,Next, we consider the operator $V_{c}$ defined as follows (for any $c\in\mathbb{R}$):
\begin{equation}
V_{c}:=\frac{\partial_{c}-c}{i},\label{Eq:inftesV}
\end{equation}
where $i:=\sqrt{-1}$ (here and throughout this appendix).\,Then, according to Equation (\ref{Eq:part}), $V_{c}$ is an unbounded self-adjoint linear operator on $\mathbb{H}_{c}$  whose domain is the same as the domain of $\partial_{c}$ (see Equation (\ref{Eq:acf})); i.e., $D(V_{c})=D(\partial_{c})$.\\

\hspace*{3mm}As a result, we obtain the following equivalent restatement of Theorem \ref{Thm:part}.

\begin{corollary}\label{CorTh}
Let $c\in \mathbb{R}$.\,Then $\partial_{c}$ is a normal operator given by
\begin{equation}
\partial_{c}=c+iV_{c}=Re(\partial_{c})+iIm(\partial_{c}),\label{Eq:partV}
\end{equation}
where $Re(\partial_{c})=c$ and $Im(\partial_{c})=V_{c}$ denote respectively the real and imaginary parts of $\partial_{c}$.\footnote{For notational simplicity, we write $c$ instead of $c$ times the identity operator of $D(\partial_{c})=D(V_{c})$.} $($Of course, it follows that $V_{c}$ is a self-adjoint operator.$)$
\end{corollary}

Our next result follows from Equation $($\ref{Eq:partV}$)$.
\begin{theorem}\label{Thm:mom}
Let $c\in \mathbb{R}$.\,Then $\frac{1}{i}\partial_{c}$ is self-adjoint if and only if $c=0$.
\end{theorem}

\hspace*{3mm}Note that the case where $c=0$ then corresponds to the usual situation of a quantum mechanical particle constrained to move on the real line $\mathbb{R}$.\, In other words, $V_{0}=\frac{1}{i}\partial_{0}$ is the standard \emph{momentum operator} acting on $\mathbb{H}_{0}=L^{2}(\mathbb{R})$.\\

\hspace*{3mm}With the notation of Corollary \ref{CorTh}, we obtain the following characterization of the spectrum of the self-adjoint `$c$-\emph{momentum operator}' $V_{c}$:
\begin{theorem}\label{Thm:SpVc}
For any $c\in \mathbb{R}$, the spectrum $\sigma(V_{c})$ of the unbounded self-adjoint operator $V_{c}=Im(\partial_{c})$ is given by
\begin{equation}
\sigma(V_{c})=\sigma_{e}(V_{c})=\mathbb{R},\label{Eq:spvc}
\end{equation}
where $\sigma_{e}(V_{c})$ \emph{denotes the essential spectrum} of $V_{c}$.
\end{theorem}

\hspace*{3mm}In particular, note that for any value of the parameter $c\in\mathbb{R}$, the spectrum of the operator $V_{c}$ coincides with the spectrum of the classic momentum operator $V_{0}$.\,In fact, we will show below that $V_{c}$ is unitarily equivalent to $V_{0}$, which is a much stronger and more precise statement.\,It follows that the point spectrum of $V_{c}$ is empty (i.e., $V_{c}$ does not have any eigenvalues) and therefore, $\sigma_{ap}(V_{c})$, the approximate point spectrum of $V_{c}$, coincides with $\sigma(V_{c})$.\,Hence, $\sigma_{ap}(V_{c})$ is also given by the right-hand side of Equation (\ref{Eq:spvc}).\,(See footnote (20) for the definition of $\sigma_{ap}$.)\\

\hspace*{3mm}Next, following \cite{HerLa1}, we outline the main steps of the proof of Corollary \ref{CorTh} and Theorem \ref{Thm:SpVc} (and hence, equivalently, of Theorems \ref{Thm:part} and \ref{Thm:spp}).\footnote{An alternative (or \textquotedblleft direct\textquotedblright) proof of Theorems \ref{Thm:part} and  \ref{Thm:spp}, not simply using the known results about the spectrum of $V_{0}$ (based on the properties of the Fourier transform and the multiplication form of the spectral theorem for self-adjoint operators), is also given in \cite{HerLa1}.}\,It is well known $($see, e.g., \cite{Sc} or vol.\,II of \cite{ReSi}$)$ that the standard momentum operator
\begin{equation}
V_{0}=\frac{1}{i}\partial_{0}=\frac{1}{i}\frac{d}{dt}
\end{equation}
is an unbounded self-adjoint operator in $L^{2}(\mathbb{R})$ since, via the Fourier transform, it is unitarily equivalent to the multiplication operator by the variable $t$ in $L^{2}(\mathbb{R})=L^{2}(\mathbb{R},dt)=\mathbb{H}_{0}$.\,Moreover, $\sigma(V_{0})=\mathbb{R}$ since by the multiplication form of the spectral theorem for unbounded self-adjoint (or, more generally, normal) operators (see, e.g., \cite{ReSi, Sc, JoLa, Ru}), $\sigma(V_{0})$ is equal to the (essential) range of the identity map $t\mapsto t$ ($t\in \mathbb{R}$), which is $\mathbb{R}$.\\
 
\hspace*{3mm}As a result, Theorem \ref{Thm:part} $($or equivalently, Corollary \ref{CorTh}$)$ can be proved by merely showing that $V_{c}=\frac{\partial_{c}-c}{i}$ is unitarily equivalent to $V_{0}$.\,More specifically, it is shown in \cite{HerLa1} (and follows from the definition of $\partial_{c}$ and of its domain, along with Equation (\ref{Eq:inftesV})) that 
\begin{equation}
V_{0}=WV_{c}W^{-1}
\end{equation}
or equivalently, 
\begin{equation}
V_{c}=W^{-1}V_{0}W,
\end{equation}
where $W:\mathbb{H}_{c}\to \mathbb{H}_{0}$ is the unitary map from $\mathbb{H}_{c}$ onto $\mathbb{H}_{0}$ defined by 
\begin{equation}
(Wf)(t)=e^{-ct}f(t)
\end{equation} 
for $f\in\mathbb{H}_{c}$, so that 
\begin{equation}
(W^{-1}g)(t)=e^{ct}g(t)
\end{equation} 
for $g\in \mathbb{H}_{0}$.\\

\hspace*{3mm}Finally, we note that in light of the above proof, for any $c\in\mathbb{R}$, $V_{c}$ is self-adjoint with spectrum $\sigma(V_{c})=\mathbb{R}$.\,(Indeed, as was recalled above, $\sigma(V_{0})=\mathbb{R}$.\,Moreover, unitary equivalence preserves the spectrum, so that $\sigma(V_{c})=\sigma(V_{0})=\mathbb{R}$.)\,Therefore, since $V_{c}=\frac{\partial_{c}-c}{i}$, we deduce that $\partial_{c}=c+iV_{c}$ is a normal unbounded operator with spectrum
\begin{equation}
\sigma(\partial_{c})=c+i\sigma(V_{c})=c+i\mathbb{R}.
\end{equation}
This establishes both Corollary \ref{CorTh} (or equivalently, Theorem \ref{Thm:part} ) and Theorem \ref{Thm:SpVc} (or equivalently, Theorem \ref{Thm:spp}).\footnote{As a result, we deduce that Theorems \ref{Thm:part} and \ref{Thm:spp} are valid without change for any $c\in \mathbb{R}$ rather for any $c\geq0$.}

\bibliographystyle{amsalpha}

\end{document}